\documentclass[preprint,3p]{elsarticle}

\usepackage{lipsum}
\makeatletter
\def\ps@pprintTitle{%
 \let\@oddhead\@empty
 \let\@evenhead\@empty
 \def\@oddfoot{}%
 \let\@evenfoot\@oddfoot}
\makeatother

\usepackage{amssymb}
\usepackage{amsthm}
\usepackage{amsbsy}
\usepackage{amsmath}
\usepackage{algorithm,algpseudocode}

\newtheorem{rem}{Remark}

\usepackage{lineno,hyperref}
\modulolinenumbers[5]

\journal{}

\begin{document}

\begin{frontmatter}
 
\title{A comprehensive deep learning-based approach to reduced order modeling of nonlinear time-dependent parametrized PDEs}

\author[1]{Stefania Fresca}
\ead{stefania.fresca@polimi.com}
\author[1]{Luca Ded\'e}
\ead{luca.dede@polimi.com}
\author[1]{Andrea Manzoni\corref{cor1}}
\ead{andrea1.manzoni@polimi.com}
\cortext[cor1]{Corresponding author}

\address[1]{MOX - Dipartimento di Matematica, Politecnico di Milano, P.zza Leonardo da Vinci 32, 20133 Milano, Italy}

\begin{abstract}
Traditional reduced order modeling techniques such as the reduced basis (RB) method (relying, e.g., on proper orthogonal decomposition (POD)) suffer from severe limitations when dealing with nonlinear time-dependent parametrized PDEs,  because of the fundamental assumption of linear superimposition of modes they are based on. For this reason, in the case of problems featuring coherent structures that propagate over time such as transport,  wave, or convection-dominated phenomena, the RB method usually yields inefficient reduced order models (ROMs) if one aims at obtaining reduced order approximations sufficiently accurate compared to the high-fidelity, full order model (FOM) solution. To overcome these limitations, in this work, we propose a new nonlinear approach to set reduced order models by exploiting deep learning (DL) algorithms. In the resulting nonlinear ROM, which we refer to as DL-ROM, both the nonlinear trial manifold (corresponding to the set of basis functions in a linear ROM) as well as the nonlinear reduced dynamics (corresponding to the projection stage in a linear ROM) are learned in a non-intrusive way by relying on DL algorithms; the latter are  trained on a set of FOM solutions obtained for different parameter values. In this paper, we show how to construct a DL-ROM for both linear and nonlinear time-dependent parametrized PDEs; moreover, we assess its accuracy on test cases   featuring different parametrized PDE problems. Numerical results indicate that DL-ROMs whose dimension is equal to the intrinsic dimensionality of the PDE solutions manifold are able to approximate the solution of parametrized PDEs  in situations where a huge number of POD modes would be necessary to achieve the same degree of accuracy.
\end{abstract}


\end{frontmatter}


\section{Introduction}
\label{sec:1}

The solution of a parametrized system of partial differential equations (PDEs) by means of a \textit{full-order model} (FOM), whenever dealing with real-time or multi-query scenarios, entails prohibitive computational costs if the FOM is high-dimensional. In the former case, the FOM solution must be computed in a very limited amount of time; in the latter one, the FOM must be solved for a huge number of parameter instances sampled from the parameter space.
Reduced order modeling techniques aim at replacing the FOM by a reduced order model (ROM), featuring a much lower dimension, still able to {express the} physical features of the problem described by the FOM. The basic assumption underlying the construction of such a ROM is that the solution of a parametrized PDE, belonging a priori to a high-dimensional (discrete) space, lies on a low-dimensional manifold embedded in this space. The goal of a ROM is then to approximate the \textit{solution manifold} -- that is, the set of all PDE solutions when the parameters vary in the parameter space -- through a suitable, approximated {\em trial manifold}.

A widespread family of reduced order modeling techniques relies on the assumption that the reduced-order approximation can be expressed by a linear combination of basis functions, built starting from a set of FOM solutions, called snapshots. Among these techniques, proper orthogonal decomposition (POD) -- equivalent to principal component analysis in statistics \cite{hastie2001theelements}, or Karhunen-Lo\`eve expansion in stochastic applications -- exploits the singular value decomposition of a suitable snapshot matrix (or the eigen-decomposition of the corresponding snapshot correlation matrix), thus yielding {\em linear} ROMs, in which the ROM approximation is given by the linear superimposition of POD modes. In this case, the solution manifold is approximated through a \textit{linear} trial manifold, that is, the ROM approximation is sought in a low-dimensional linear trial subspace.

Projection-based methods are linear ROMs in which the ROM approximation of the PDE solution, for any new parameter value, results from the solution of a low-dimensional (nonlinear, dynamical) system, whose unknowns are the ROM degrees of freedom (or generalized coordinates). Despite the PDE (and thus the FOM) {being} linear or not, the operators appearing in the ROM are obtained by imposing that the projection of the FOM residual evaluated on the ROM trial solution is orthogonal to a low-dimensional, linear test subspace, which might coincide with the trial subspace. Hence, no matter whether the PDE is linear or not,  {the resulting ROM is {\em linear} since the reduced dynamics is obtained through a  projection onto a linear subspace} \cite{benner2017model,benner2015asurvey,quarteroni2016reduced}. However, linear ROMs show severe computational bottlenecks when dealing with problems featuring coherent structures (possibly dependent on parameters) that propagate over time, namely in transport and wave-type phenomena, or convection-dominated flows. In these cases, the dimension of the linear trial manifold can easily become extremely large if compared to the intrinsic dimension of the solution manifold for the sake of accuracy, thus compromising the ROM efficiency. To overcome this bottleneck, {\em ad-hoc} extensions of the POD strategy have been considered, towards nonlinear approaches to build a ROM \cite{ohlberger2016reduced, pagani2018numerical}.

In this paper we propose a computational, non-intrusive approach based on deep learning (DL) algorithms to deal with the construction of efficient ROMs (which we refer to as DL-ROMs) in order to tackle parameter-dependent PDEs; in particular, we consider PDEs that feature wave-type phenomena. A comprehensive framework is presented for the global approximation of the map $(t, \boldsymbol{\mu}) \mapsto {\mathbf{u}}_h(t, \boldsymbol{\mu})$, where $t \in (0,T)$ denotes time, $\boldsymbol{\mu} \in \mathcal{P} \subset \mathbb{R}^{n_{\mu}}$ a vector of input parameters and $ {\mathbf{u}}_h(t, \boldsymbol{\mu}) \in \mathbb{R}^{N_h}$ the solution of a large-scale dynamical system arising from the space discretization of a parametrized, time-dependent (non)linear PDE. Several recent works have shown possible applications of DL techniques to parametrized PDEs -- thanks to their approximation capabilities, their extremely favorable computational performances during online testing phases, and their relative easiness of implementation --  both from a theoretical \cite{kutyniok2019atheoretical} and a computational standpoint. 
Regarding this latter aspect, artificial neural networks (ANN), such as feedforward neural networks, 
have been employed to model the reduced dynamics in a data-driven \cite{regazzoni2019machinelearningfor}, and less intrusive way (avoiding, e.g., the costs entailed by projection-based ROMs), but still relying on a linear trial manifold built, e.g., through POD. For instance, in \cite{guo2018reduced, guo2019data, hestaven2018non-intrusive, san2018neural} the solution of a (nonlinear, time-dependent) ROM for any new parameter value has been replaced by the evaluation of ANN-based regression models; similar ideas can be found, e.g., in \cite{kani2017dr-rnn,mohan2018adeep,wan2018data}. Few attempts  have been made in order to describe the reduced trial manifold where the approximation is sought (avoiding, e.g., the linear superimposition of POD modes) through ANNs, see, e.g., \cite{gonzalez2018deep, carlberg2018model}.

For instance, a projection-based ROM technique has been introduced in \citep{carlberg2018model}, in which the FOM system is projected onto a nonlinear trial manifold identified by means of the decoder function of a convolutional autoencoder neural network. However, the ROM is derived by minimizing a residual formulation, for which the quasi-Newton method herein employed requires the computation of an approximated Jacobian of the residual at each time step. A ROM technique based on a deep convolutional recurrent autoencoder has been proposed in \cite{gonzalez2018deep}, where a reduced  trial manifold is obtained by means of a convolutional autoencoder; the latter is then used to train a Long Short-Term Memory (LSTM) neural network modeling the reduced dynamics. However, no explicit parameter dependence in the PDE problem is considered, apart from $\boldsymbol{\mu}$-dependent initial data, and the LSTM is trained on reduced approximations obtained through the encoder function of the autoencoder. Another promising application of machine and deep learning techniques within a ROM framework deals with the efficient evaluation of reduced error models, see, e.g., \cite{ freno2018machine, paganicarlbergmanzoni2019, parish2019time, trehan2017error}. \\ 

Our goal is to set up nonlinear ROMs whose dimension is nearly equal (if not equal) to  the intrinsic dimension of the solution manifold that we aim at approximating. Our DL-ROM approach  combines and  improves the techniques introduced in \cite{gonzalez2018deep,carlberg2018model} by shaping an all-inclusive DL-based ROM technique, where we  both {\em (i)} construct the reduced trial manifold and {\em (ii)} model the reduced dynamics on it employing ANNs. The former task is achieved by using  the decoder function of a convolutional autoencoder; the latter task is instead carried out by considering a feedforward neural network {and the encoder function of a convolutional autoencoder}. Moreover, we set up a computational procedure performing the training of both network architectures simultaneously, by minimizing a loss function that weights  two terms, one dedicated to  each single task.  In this respect, we are able to design a flexible framework capable to handle parameters affecting both PDE operators and data, which avoids both the expensive projection stage of \cite{carlberg2018model} and the training of a more expensive LSTM network.
In our technique, the intrusive construction of a ROM is replaced by the evaluation of the ROM generalized coordinates through a deep feedforward neural network taking  only $(t, \boldsymbol{\mu})$ as inputs. The proposed technique is purely data-driven, that is, it only relies on the computation of a set of FOM snapshots -- in this respect, DL does not replace the high-fidelity FOM as, e.g., in the works by Karniadakis and coauthors \cite{raissi2018hidden, raissi2017physics1, raissi2017physics2, raissi2019physics, raissi2018deep}; rather, DL techniques are built upon it, to enhance the repeated evaluation of the FOM for different values of the parameters. 

The structure of the paper is as follows. In \autoref{sec:2} we show how to generate nonlinear ROMs by reinterpreting the classical ideas behind linear ROMs for parametrized PDEs. In \autoref{sec:3} we detail the construction of the proposed DL-ROM, whose accuracy is numerically assessed in \autoref{sec:4} by considering three different test cases of increasing complexity (both with respect to the parametric dependence and the nature of the PDE). Finally, some conclusions are drawn in \autoref{sec:5}. A quick overview of useful facts about deep feedforward, convolutional  and autoencoders neural networks is reported in \ref{sec:A} to make the paper self-contained.

\section{From linear to nonlinear dimensionality reduction}
\label{sec:2}

Starting from the well-known setting of linear (projection-based) ROMs, in this section we  generalize this task to the case of  nonlinear ROMs.

\subsection{Problem formulation}

We formulate the construction of ROMs in algebraic terms, starting from the high-fidelity (spatial) approximation of nonlinear, time-dependent, parametrized PDEs. By introducing suitable space discretizations techniques (such as, e.g., the Finite Element Method, Isogeometric Analysis  or the Spectral Element Method) the high-fidelity, full order model (FOM) can be expressed as a nonlinear parametrized dynamical system. Given $\boldsymbol{\mu} \in \mathcal{P}$, we aim at solving the initial value problem
\begin{equation}
\label{FOM}
\begin{cases}
\mathbf{\dot{u}}_h(t;\boldsymbol{\mu}) = \mathbf{f}(t, \mathbf{u}_h(t;\boldsymbol{\mu}); \boldsymbol{\mu}) \quad t \in (0, T)\\
\mathbf{u}_h(0;\boldsymbol{\mu})=\mathbf{u}_0(\boldsymbol{\mu}),
\end{cases}
\end{equation}
where  the parameter space $\mathcal{P} \subset \mathbb{R}^{n_{\boldsymbol{\mu}}}$ is a bounded and closed set,  $\mathbf{u}_h:[0,T) \times \mathcal{P} \rightarrow \mathbb{R}^{N_h}$ is the parametrized  solution of (\ref{FOM}),  $\mathbf{u}_0 : \mathcal{P} \rightarrow \mathbb{R}^{N_h}$ is the initial datum and $\mathbf{f} : (0,T) \times  \mathbb{R}^{N_h}  \times \mathcal{P} \rightarrow \mathbb{R}^{N_h}$ is a (nonlinear) function, encoding the system dynamics. The FOM dimension $N_h$  is related with the finite dimensional subspaces introduced for the space discretization of the PDE -- here $h>0$ usually denotes a discretization parameter, such as the maximum diameter of elements in a computational mesh --  and can be extremely small whenever the PDE problem shows complex physical behaviors and/or high degrees of accuracy are required to its solution. The parameter $\boldsymbol{\mu} \in \mathcal{P}$ may represent physical or geometrical properties of the system, like, e.g., material properties, initial and boundary conditions, or the shape of the domain.  In order to solve problem (\ref{FOM}), suitable time discretizations are employed, such as backward differentiation formulas \cite{quarteroni2008matematica}.

\begin{figure}[t!]
\centering
\includegraphics[width=0.475\textwidth]{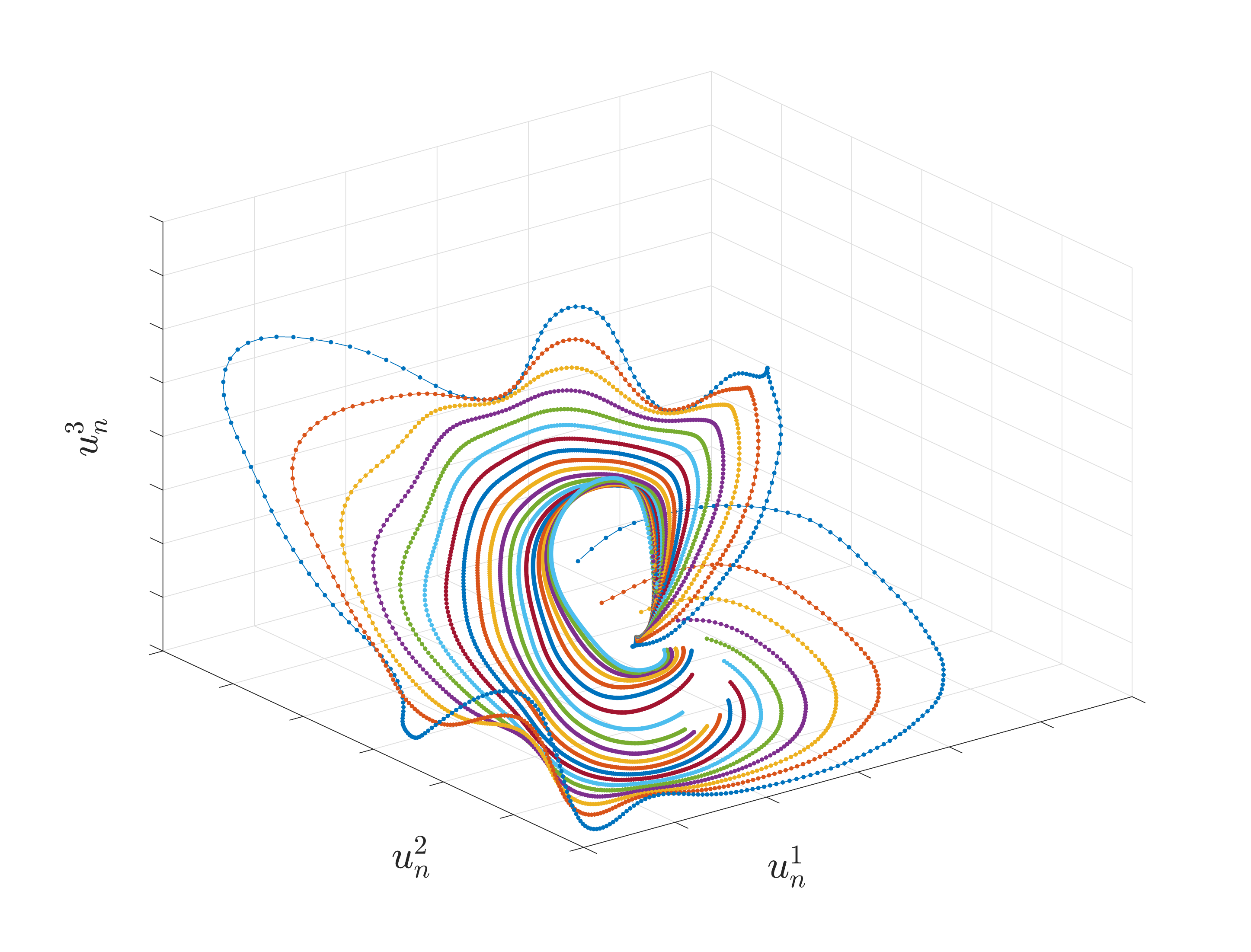}
\caption{Example of a two-dimensional manifold embedded in $\mathbb{R}^3$. Each curve represents the time-evolution of the first three components of the solution of a (nonlinear) parametrized PDE for a fixed parameter value $\boldsymbol{\mu}$. }
\label{manifold}
\end{figure}

Our goal is the efficient numerical approximation of the whole set
\begin{equation}
\mathcal{S}_h = \{ \mathbf{u}_h(t ; \boldsymbol{\mu} ) \; | \; t \in [0, T) \; \textnormal{and} \; \boldsymbol{\mu} \in \mathcal{P} \subset \mathbb{R}^{n_{\mu}} \} \subset \mathbb{R}^{N_h},
\label{solution_manifold}
\end{equation}
of solutions to problem (\ref{FOM}) when $(t ; \boldsymbol{\mu} )$ varies in $[0, T) \times \mathcal{P} $, also referred to as {\em solution manifold}  (a sketch is provided  in \figurename~\ref{manifold}). 
Assuming that, for any given parameter $\boldsymbol{\mu} \in \mathcal{P}$, problem (\ref{FOM}) admits a unique solution,  for each $t \in (0,T)$, the intrinsic dimension of the solution manifold is at most $n_{\boldsymbol{\mu}} + 1  \ll N_h$, where $n_{\boldsymbol \mu}$ is the number of parameters (time plays the role of an additional coordinate). This means that each point $\mathbf{u}_h(t ; \boldsymbol{\mu} )$ belonging to $\mathcal{S}_h$  is completely defined in terms of at most $n_{\boldsymbol{\mu}} + 1$ intrinsic coordinates, or equivalently, the tangent space to the manifold at any given $\mathbf{u}_h(t ; \boldsymbol{\mu} )$ is spanned by $n_{\boldsymbol{\mu}} + 1$ basis vectors.

\subsection{Linear dimensionality reduction: projection-based ROMs}

The most common way to build a ROM for the efficient approximation of problem \eqref{FOM} relies on the introduction of a \textit{reduced linear  trial manifold}, that is of a subspace $\tilde{\mathcal{S}}_n = \textnormal{Col}(V)$ of dimension $n \ll N_h$, spanned by the $n$ columns of a matrix $V \in \mathbb{R}^{N_h \times n}$. Hence, a linear ROM looks for an approximation $\mathbf{\tilde{u}}_h(t; \boldsymbol \mu) \approx \mathbf{u}_h(t; \boldsymbol \mu)$ in the form
\begin{equation}
\label{linear_reconstructed_solution}
\mathbf{\tilde{u}}_h(t; \boldsymbol \mu) = V\mathbf{u}_n(t; \boldsymbol{\mu}) , 
\end{equation}
where 
$\mathbf{\tilde{u}}_h : [0,T) \times \mathcal{P} \rightarrow \tilde{\mathcal{S}}_n$. 
Here $\mathbf{u}_n(t; \boldsymbol{\mu}) \in \mathbb{R}^n$ for each $t \in [0,T)$, $\boldsymbol{\mu} \in \mathcal{P}$ denotes the vector of intrinsic coordinates (or degrees of freedom) of the ROM approximation;  note that the map
\[
\boldsymbol{\Psi}_{h} : \mathbb{R}^n 	\rightarrow \mathbb{R}^{N_h}, \qquad
 \mathbf{s}_n  \mapsto \tilde{\mathbf{s}}_h = V \mathbf{s}_n
 \]
that, given the (low-dimensional) intrinsic coordinates, returns the (high-dimensional) approximation of the FOM solution ${\mathbf{u}}_h(t; \boldsymbol{\mu})$, is linear.

Proper Orthogonal Decomposition (POD) is one of the most widely employed techniques to generate the linear trial manifold \cite{quarteroni2016reduced}. Considering a set of $N_{train}$ instances of the parameter $\boldsymbol \mu \in \mathcal{P}$, we introduce the snapshot matrix $S \in \mathbb{R}^{N_h \times N_s}$ defined as
\begin{equation}
\label{snapshot_matrix}
S = \left[ \mathbf{u}(t^1; \boldsymbol \mu_1) \; | \; \ldots \; | \; \mathbf{u}(t^{N_t}; \boldsymbol \mu_1) \; | \; \ldots \; | \;
\mathbf{u}(t^1 ; \boldsymbol \mu_{N_{train}}) \; | \; \ldots \; | \; \mathbf{u}(t^{N_t} ; \boldsymbol \mu_{N_{train}}) \right],
\end{equation}
where we have introduced a partition of the time interval $[0,T]$ in $N_t$ time steps $\{t^k\}_{{k=1}}^{N_t}$, $t^k = k \Delta t$, of size $\Delta t =T/N_t$  {and $N_s = N_{train}N_t$}. Moreover, let us introduce  a symmetric and positive definite matrix $X_h \in \mathbb{R}^{N_h \times N_h}$ encoding a suitable norm (e.g., the energy norm) on the high-dimensional space and  admitting a Cholesky factorization $X_h = H^T H$. POD computes the Singular Value Decomposition (SVD) of   $H  S$,
\begin{equation*}
H S = U \Sigma Z^T,
\end{equation*}
where $U = [\boldsymbol{\zeta}_1| \ldots | \boldsymbol{\zeta}_{N_h}] \in \mathbb{R}^{N_h \times N_h}$, $Z = [\boldsymbol{\psi}_1| \ldots | \boldsymbol{\psi}_{N_s}] \in \mathbb{R}^{N_s \times N_s}$ and $\Sigma = \textnormal{diag}(\sigma_1, \ldots, \sigma_r) \in \mathbb{R}^{N_h \times N_s}$ with $\sigma_1 \geq \sigma_2 \geq \ldots \geq \sigma_r$, and $r \leq \min(N_h, N_s)$, and sets the columns of $V$ in terms of the first $n$ left singular vectors of $S$ that is,  $V = [H^{-1} \boldsymbol{\zeta}_1 | \ldots | H^{-1}  \boldsymbol{\zeta}_n]$.
By construction, the columns of $V$ are orthonormal (with respect to the scalar product $( \, \cdot \, , \, \cdot \,)_{X_h}$) and
among all possible $n$-dimensional subspaces spanned by the column of a matrix $W \in \mathbb{R}^{N_h 	\times n}$, $V$ provides the best reconstruction of the snapshots, that is,
\begin{equation}
\label{minimization_problem_POD}
  \sum_{i=1}^{N_{train}} \sum_{k=1}^{N_t} \| \mathbf{u}(t^k;\boldsymbol \mu_i) - VV^T X_h \mathbf{u}(t^k;\boldsymbol \mu_i) \|_{X_h}^2 = \min_{W \in \mathcal{V}_{n}} \sum_{i=1}^{N_{train}}
\sum_{k=1}^{N_t} \| \mathbf{u}(t^k; \boldsymbol \mu_i) - WW^T X_h \mathbf{u}(t^k; \boldsymbol \mu_i) \|_{X_h}^2,
\end{equation}
where $\mathcal{V}_{n} = \{ W \in \mathbb{R}^{N_h \times n} : W^T X_h W = I \}$. For this reason, we refer to $V V^T X_h \mathbf{u}_h(t; \boldsymbol{\mu})$ as to the optimal-{POD} reconstruction of $\mathbf{u}_h(t; \boldsymbol{\mu})$ onto a reduced subspace of dimension $n < N_h$.

In order to model the reduced dynamics of the system, that is, the time-evolution of the generalized coordinates $\mathbf{u}_n(t; \boldsymbol{\mu})$, we can replace $\mathbf{u}_h(t; \boldsymbol{\mu})$ by (\ref{linear_reconstructed_solution}) in system (\ref{FOM}),
\begin{equation}
\label{FOM_linear_reconstruction}
\begin{cases}
V \displaystyle \mathbf{\dot{u}}_n(t; \boldsymbol{\mu})  = \mathbf{f}(t, V\mathbf{u}_n(t; \boldsymbol{\mu});  \boldsymbol{\mu}) \quad t \in (0, T) \\
V \mathbf{u}_n(0; \boldsymbol \mu) = \mathbf{u}_0(\boldsymbol \mu),
\end{cases}
\end{equation}
and impose that the residual
\begin{equation}
\label{residual}
\mathbf{r}_h(V\mathbf{u}_n (t; \boldsymbol{\mu})) = V \mathbf{\dot{u}}_n(t; \boldsymbol{\mu}) - \mathbf{f}(t, V\mathbf{u}_n (t; \boldsymbol{\mu}); \boldsymbol{\mu})
\end{equation}
associated to the first equation of (\ref{FOM_linear_reconstruction}) is  orthogonal to a $n$-dimensional subspace spanned by the column of a matrix $Y \in \mathbb{R}^{N_h \times n}$, that is, $Y^T \mathbf{r}_h(V\mathbf{u}_n) = {\bf 0}$. This condition yields the following ROM
\begin{equation}
\label{ROM}
\begin{cases}
Y^T V \mathbf{\dot{u}}_n(t;\boldsymbol{\mu}) = Y^T \mathbf{f}(t, V\mathbf{u}_n(t;\boldsymbol{\mu}); \boldsymbol{\mu}) \quad t \in (0, T)\\
\mathbf{u}_n(0;\boldsymbol{\mu})= (Y^T V)^{-1} Y^T\mathbf{u}_0(\boldsymbol{\mu}).
\end{cases}
\end{equation}
In the case $Y=V$ a Galerkin projection is performed, while the case $Y \neq V$  yields a more general Petrov-Galerkin projection. Note that choosing $Y$ such that $Y^T V = I \in \mathbb{R}^{N_h \times N_h}$ does not automatically ensure ROM stability on long time intervals.

The RB method under the form of either Galerkin-POD or Petrov-Galerkin-POD methods has been successfully applied to a broad range of parametrized time-dependent (non)linear problems (see, e.g., \citep{pagani2018numerical,manzoni2016accurate}) however it provides low-dimensional subspaces of dimension $n \gg n_{\mu} + 1$ much larger than the intrinsic dimension of the solution manifold --  relying on a linear, global trial manifold thus represent a major bottleneck to computational efficiency \citep{ohlberger2016reduced,pagani2018numerical}. This is the case, for instance,  of hyperbolic problems, for which the RB method is not able in practice to significantly decrease the dimensionality of the problem. The same difficulty might also affect the use of hyper-reduction techniques, such as the (discrete) empirical interpolation \cite{barrault2004anempirical, chaturantabut2010nonlinear}, mandatory in order to assemble the operators appearing in the ROM \eqref{ROM} without relying on expensive $N_h$-dimensional arrays. See, e.g., \cite{FGMQ_19} for further details.

\subsection{Nonlinear dimensionality reduction}
\label{sec:nonlinear_dimred}

A first attempt to overcome the computational bottleneck  entailed by the use of a linear, global trial manifold is to build a {\em piecewise} linear trial manifold, using local reduced bases whose dimension is smaller than the one of the global linear trial manifold. Clustering algorithms applied on a set of snapshots can be employed to partition them into $N_c$ clusters from which POD can extract a subspace of reduced dimension; the ROM is then obtained by following the strategy described above on each cluster separately, see, e.g. \cite{amsallem2012nonlinear,Amsallem2015}. An alternative approach based on classification binary trees has been introduced in \cite{Amsallem2016}. These strategies have been employed (and compared) in \cite{pagani2018numerical} in order to solve parametrized problems in cardiac electrophysiology. Using a piecewise linear trial manifold partially overcomes the limitation of a linear dimensionality reduction technique as POD, yet employing local bases of dimension much higher than the intrinsic dimension of the solution manifold $\mathcal{S}_h$. 
An approach based on a dictionary of solutions, computed offline, has been developed in \cite{abgrall_amsallem} as an alternative to using a truncated reduced basis based on POD, together with an online $L^1$-norm minimization of the residual.

Other possible options involving nonlinear transformations of modes might rely on a reconstruction of the POD modes at each time step using Lax pairs \cite{GerbeauLombardi2014}, on the solution of Monge-Kantorovich optimal transport problems \cite{PhysRevE.89.022923}, on a problem-dependent change of coordinates requiring the solution of an optimization problem repeatedly \cite{cagniart2019model}, on shifted POD modes \cite{reiss2018shifted} after multiple transport velocities have been identified and separated, or again basis updates are derived from querying the full model at a few selected spatial coordinates \cite{peherstorfer2018}. Despite providing remarkable improvements compared to the {\em classic}  {(Petrov-)Galerkin-}POD approach, all these strategies exhibit some drawbacks, such as: {\em (i)} the high computational costs entailed during the online testing evaluation stage of the ROM -- which is not restricted to the intensive offline training stage; {\em (ii)} performances and settings are  highly dependent on the problem at hand; {\em (iii)} the need to deal only with a linear superimposition of modes (which characterizes linear ROMs), yielding  low-dimensional spaces whose dimension is still (much) higher than the intrinsic dimension of the solution manifold. 

Motivated by the need of avoiding the drawbacks of linear ROMs and setting a general paradigm for the construction of efficient, extremely low-dimensional ROMs, we resort to nonlinear dimensionality reduction techniques. Similarly to \cite{gonzalez2018deep, carlberg2018model}, we build a nonlinear ROM to approximate $ \mathbf{u}_h(t; \boldsymbol \mu) \approx \mathbf{\tilde{u}}_h(t; \boldsymbol \mu)$ by 
\begin{equation}
\mathbf{\tilde{u}}_h(t; \boldsymbol \mu) = \boldsymbol{\Psi}_h(\mathbf{u}_n(t; \boldsymbol{\mu})), 
\label{reconstructed_solution}
\end{equation}
where  $\boldsymbol{\Psi}_h : \mathbb{R}^{n} \rightarrow \mathbb{R}^{N_h}$, $\boldsymbol{\Psi}_h : \mathbf{s}_n \mapsto  \boldsymbol{\Psi}_h(\mathbf{s}_n)$,  $n \ll N_h$, is a nonlinear, differentiable function. As a matter of fact, 
the solution manifold $\mathcal{S}_h$ is approximated by a \textit{reduced nonlinear trial manifold}
\begin{equation}
\label{g}
\tilde{\mathcal{S}}_n = \{ \boldsymbol{\Psi}_h(\mathbf{u}_n(t; \boldsymbol{\mu}))
\; | \; \mathbf{u}_n(t; \boldsymbol{\mu}) \in \mathbb{R}^{n}, \ t \in [0, T) \; \textnormal{and} \; \boldsymbol{\mu} \in \mathcal{P} \subset \mathbb{R}^{n_{\mu}} \} \subset \mathbb{R}^{N_h}
\end{equation}
so that  $\mathbf{\tilde{u}}_h : [0,T) \times \mathcal{P} \rightarrow \tilde{\mathcal{S}}_n$. As before, $\mathbf{u}_n: [0,T) \times \mathcal{P} \rightarrow \mathbb{R}^{n}$ denotes the vector-valued function of two arguments representing the  intrinsic coordinates of the ROM approximation.
Our goal is to set a ROM whose dimension $n$ is as close as possible to the   intrinsic dimension $n_{\boldsymbol{\mu}} + 1$ of the solution manifold $\mathcal{S}_h$, i.e. $n \geq n_{\boldsymbol{\mu}} + 1$, in order to correctly capture the solution of the dynamical system by containing the size of  the approximation spaces \citep{carlberg2018model}.

To model the relationship between each couple $(t, \boldsymbol{\mu}) \mapsto \mathbf{u}_n(t, \boldsymbol{\mu})$, and  to describe the system dynamics on the reduced nonlinear trial manifold $\tilde{\mathcal{S}}_n$ in terms of the intrinsic coordinates, we consider a nonlinear map  under the form
\begin{equation}
\label{h}
\mathbf{u}_n(t; \boldsymbol \mu) = \boldsymbol{\Phi}_n(t, \boldsymbol \mu), 
\end{equation}
where $\boldsymbol{\Phi}_n : [0, T) \times \mathbb{R}^{n_{\boldsymbol \mu} + 1} \rightarrow \mathbb{R}^{n}$ is a differentiable nonlinear function. No additional assumptions such as, e.g., the (exact, or approximate) affine $\boldsymbol{\mu}$-dependence as in the RB method, are needed.

\section{A deep learning-based reduced order model (DL-ROM)}
\label{sec:3}

We now detail the construction of the proposed nonlinear ROM. In this respect, we define the functions $\boldsymbol{\Psi}_h$ and $\boldsymbol{\Phi}_n$ in (\ref{reconstructed_solution}) and (\ref{h}) by means of deep learning (DL) algorithms, exploiting neural network architectures. This choice is motivated by their ability of effectively approximating nonlinear maps, and by their ability to learn from data and generalize to unseen data. On the other hand, DL models enable us to build non-intrusive, completely data-driven, ROMs, since their construction only requires to access the dataset, the parameter values and the snapshot matrix, but not the FOM arrays appearing in (\ref{FOM}).

\begin{figure}[ht]
\centering
\includegraphics[scale=0.45]{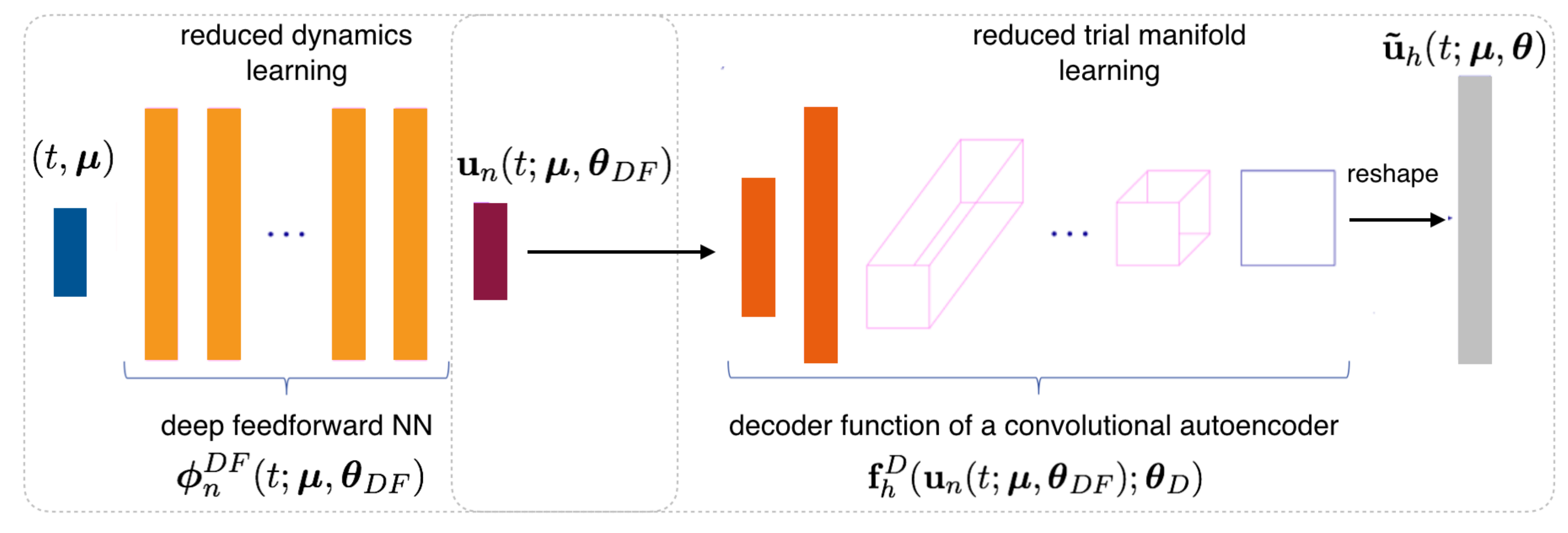}
\caption{DL-ROM architecture (online stage, testing).}
\label{architecture_DL-ROM}
\end{figure}

The DL-ROM technique that we develop in this paper is composed by two main blocks responsible, respectively, for the \textit{reduced dynamics learning} and the \textit{reduced trial manifold learning} (see \figurename~\ref{architecture_DL-ROM}). Hereon, we denote by $N_{train}$, $N_{test}$ and $N_t$ the number of training-parameter instances, of testing-parameter instances and time instances, respectively, and we set $N_{s} = N_{train} N_t$. The dimension of both the FOM solution and the ROM approximation is $N_h$, while  $n$ denotes the number of intrinsic coordinates, with $n \ll N_h$.

For the description of the system dynamics on the reduced nonlinear trial manifold (which we refer to as  {\em reduced dynamics learning}), we employ a {\em deep feedforward neural network} (DFNN) with $L$ layers, that is, we
define  the function $\boldsymbol{\Phi}_n$ in definition (\ref{h}) as 
\begin{equation}
\label{h_learning}
\boldsymbol{\Phi}_n(t; \boldsymbol \mu, \boldsymbol{\theta}_{DF}) = \boldsymbol{\phi}_n^{DF}(t; \boldsymbol \mu, \boldsymbol{\theta}_{DF}),
\end{equation}
thus yielding the map
\[
(t, \boldsymbol \mu) \mapsto \mathbf{u}_n(t; \boldsymbol \mu, \boldsymbol{\theta}_{DF}) = \boldsymbol{\phi}_n^{DF}(t; \boldsymbol \mu, \boldsymbol{\theta}_{DF}), 
\]
where $\boldsymbol{\phi}_n^{DF}$ takes the form (\ref{DFNN}), {with} $t \in [0,T]$, and results from the subsequent composition of a nonlinear activation function,  applied to a linear transformation of the input, $L$ times. Here $\boldsymbol \mu \in \mathcal{P} \subset \mathbb{R}^{n_{\boldsymbol \mu}}$ and $\boldsymbol{\theta}_{DF}$ denotes the vector of {parameters} of the DFNN.

Regarding instead the description of the reduced nonlinear trial manifold $\tilde{\mathcal{S}}_n$ defined in \eqref{g} (which we refer to as  {\em reduced trial manifold learning}), we employ the {\em decoder function of a convolutional autoencoder} (AE), that is, we define the function $\boldsymbol{\Psi}_h$ appearing in  (\ref{reconstructed_solution}) and (\ref{g}) as
\begin{equation}
\label{g_learning}
\boldsymbol{\Psi}_h(\mathbf{u}_n(t; \boldsymbol{\mu}, {\boldsymbol{\theta}_{DF}}); \boldsymbol{\theta}_{D}) = \mathbf{f}_{h}^D(\mathbf{u}_n(t; \boldsymbol{\mu}, {\boldsymbol{\theta}_{DF}}); \boldsymbol{\theta}_{D}),
\end{equation}
thus yielding the map
\[
\mathbf{u}_n(t; \boldsymbol{\mu}, {\boldsymbol{\theta}_{DF}}) \mapsto \tilde{\mathbf{u}}_h(t; \boldsymbol{\mu}, \boldsymbol{\theta} ) = \mathbf{f}_{h}^D(\mathbf{u}_n(t; \boldsymbol{\mu},  {\boldsymbol{\theta}_{DF}}); \boldsymbol{\theta}_{D}),
\]
where $\mathbf{f}_{h}^D$ results from the composition of several layers, some of which of convolutional type, overall depending on the vector $\boldsymbol{\theta}_{D}$ of parameters of the decoder function.

Combining the two former stages, the DL-ROM approximation is then given by  
\begin{equation}
\label{DL-ROM}
\mathbf{\tilde{u}}_h(t; \boldsymbol{\mu}, \boldsymbol{\theta}) = \mathbf{f}_{h}^D (\boldsymbol{\phi}_n^{DF}(t; \boldsymbol{\mu}, \boldsymbol{\theta}_{DF}); \boldsymbol{\theta}_{D}),
\end{equation}
where $\boldsymbol{\phi}_n^{DF} (\cdot; \cdot, \boldsymbol{\theta}_{DF}) : \mathbb{R}^{(n_{\boldsymbol{\mu}}+1)} \rightarrow \mathbb{R}^{n}$ and $\mathbf{f}_{h}^D (\cdot; \boldsymbol{\theta}_{D}) : \mathbb{R}^{n} \rightarrow \mathbb{R}^{N_h}$ are defined as in (\ref{h_learning}) and (\ref{g_learning}), respectively, and $\boldsymbol{\theta} = (\boldsymbol{\theta}_{DF}, \boldsymbol{\theta}_{D})$ are the parameters defining the neural network. The architecture of DL-ROM is shown in \figurename~\ref{architecture_DL-ROM}.

Computing the ROM approximation (\ref{DL-ROM}) for any new value of $\boldsymbol{\mu} \in \mathcal{P}$, at any given time, requires to evaluate the map $(t, \boldsymbol \mu) \rightarrow \mathbf{\tilde{u}}_h(t; \boldsymbol{\mu}, \boldsymbol{\theta})$  at the testing stage, once the parameters 
$\boldsymbol{\theta} = (\boldsymbol{\theta}_{DF}, \boldsymbol{\theta}_{D})$ have been determined, once and for all, during the training stage. The training stage consists in solving an  optimization problem (in the  variable $\boldsymbol{\theta}$) after a set of snapshots of the FOM have been computed. More precisely, provided the parameter matrix $M \in \mathbb{R}^{(n_{\boldsymbol \mu} + 1) \times N_{s}}$ defined as
\begin{equation}
M = [(t^1, \boldsymbol{\mu}_1) | \ldots | (t^{N_t}, \boldsymbol{\mu}_1) | \ldots | (t^1, \boldsymbol{\mu}_{N_{train}}) | \ldots | (t^{N_t}, \boldsymbol{\mu}_{N_{train}})],
\end{equation}
and the snapshot matrix $S$, defined in (\ref{snapshot_matrix}), we solve the problem:  find the optimal parameters $\boldsymbol{\theta}^*$ solution of
\begin{equation}
\label{minimization_problem_DL-ROM}
  \mathcal{J}(\boldsymbol{\theta}) =  \frac{1}{N_s}\sum_{i=1}^{N_{train}}\sum_{k=1}^{N_{t}} \mathcal{L}(t^k, \boldsymbol{\mu}_i;  {\boldsymbol{\theta}}) \rightarrow \min_{\boldsymbol{\theta}} 
\end{equation}
where
\begin{equation}
\label{loss_DL-ROM}
\mathcal{L}(t^k, \boldsymbol{\mu}_i;  {\boldsymbol{\theta}}) = \frac{1}{2}\| \mathbf{u}_h(t^k; \boldsymbol{\mu}_i) - \mathbf{\tilde{u}}_h(t^k; \boldsymbol{\mu}_i, \boldsymbol{\theta})\|^2 =
\frac{1}{2}\| \mathbf{u}_h(t^k; \boldsymbol{\mu}_i) - \mathbf{f}_{h}^D (\boldsymbol{\phi}_n^{DF}(t^k; \boldsymbol{\mu}_i, \boldsymbol{\theta}_{DF}); \boldsymbol{\theta}_{D})\|^2.
\end{equation}

To solve the optimization problem (\ref{minimization_problem_DL-ROM})-(\ref{loss_DL-ROM}) we use the ADAM algorithm \cite{kingma2015adam} which is a stochastic gradient descent method \citep{robbins1951astochastic} computing an adaptive approximation of the first and second momentum of the gradients of the loss function. In particular, it computes exponentially weighted moving averages of the gradients and of the squared gradients. We set the starting learning rate to $\eta = 10^{-4}$, the batch size to $N_b = 20$ and the maximum number of epochs to $N_{epochs} = 10000$. We perform cross-validation, in order to tune the hyper-parameters of the DL-ROM, by splitting the data in training and validation and following a proportion 8:2. Moreover, we implement an early-stopping  regularization technique  to reduce overfitting \citep{goodfellow2016deep}. In particular, we stop the training if the loss does not decrease over 500 epochs. As nonlinear activation function we employ the ELU function \citep{clevert2015fast} defined as
\begin{equation*}
\sigma(z) =
\begin{cases}
z & z \geq 0 \\
\exp(z) - 1 & z < 0.
\end{cases}
\end{equation*}
No activation function is applied at the last convolutional layer of the decoder neural network, as usually done when dealing with autoencoders. The parameters, weights and biases, are initialized through the He uniform initialization \citep{he2015delving}.

As we rely on a convolutional autoencoder to define the function $\boldsymbol{\Psi}_h$, we also exploit the encoder function 
\begin{equation}
\label{encoder}
 {\mathbf{\tilde{u}}_n}(t; \boldsymbol{\mu}, \boldsymbol{\theta}_{E}) = \mathbf{f}_{n}^E(\mathbf{u}(t; \boldsymbol{\mu}); \boldsymbol{\theta}_{E}),
\end{equation}
which maps each FOM solution associated to the pairs $(t; \boldsymbol{\mu}) \in Col(M)$ provided as inputs to the feed-forward neural network (\ref{h_learning}), onto a low-dimensional representation $\tilde{\mathbf{u}}_n(t; \boldsymbol{\mu}, \boldsymbol{\theta}_E)$ depending on the parameters vector $\boldsymbol{\theta}_E$ defining the encoder function.

\begin{figure}[ht!]
\centerline{
\includegraphics[width=0.85\textwidth]{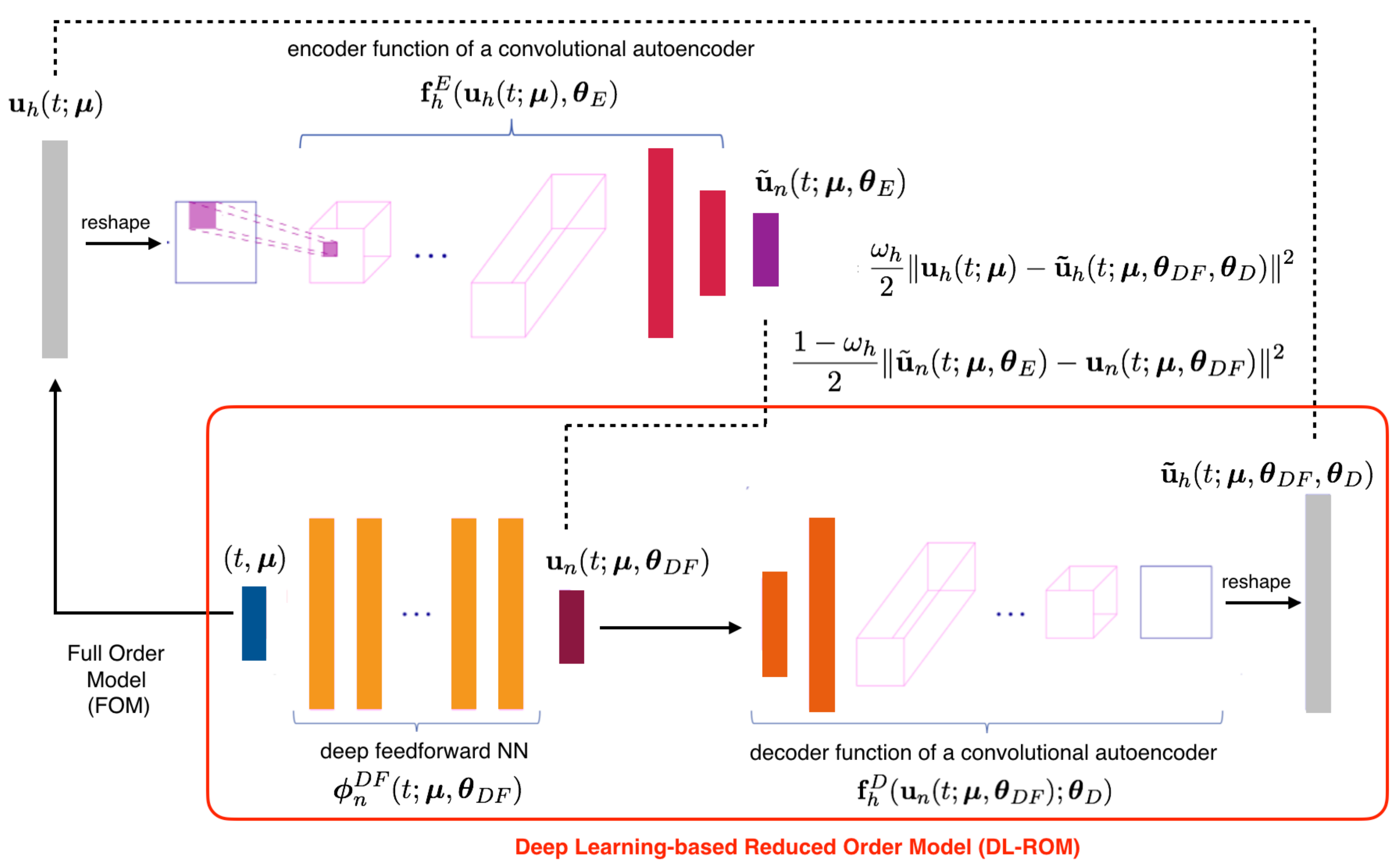}
}
\caption{DL-ROM architecture (offline stage, training and validation).}
\label{architecture_encoder}
\end{figure}

Indeed, the actual architecture of DL-ROM that is used only during the training and the validation phases, but not during testing, is the one shown in \figurename~\ref{architecture_encoder}. 
In practice, we add to the  architecture of the DL-ROM introduced above the encoder function of the convolutional autoencoder.  This {produces} an additional term in the \emph{per-example} loss function (\ref{loss_DL-ROM}), thus calling the following optimization problem to be solved:
\begin{equation}
\label{minimization_problem_encoder}
\min_{\boldsymbol{\theta}} \mathcal{J}(\boldsymbol{\theta}) = \min_{\boldsymbol{\theta}} \frac{1}{N_s}\sum_{i=1}^{N_{train}} \sum_{k=1}^{N_t} \mathcal{L}(t^k, \boldsymbol \mu_i; \boldsymbol{\theta}), 
\end{equation} 
where
\begin{equation}
\label{loss_encoder}
\mathcal{L}(t^k, \boldsymbol{\mu}_i;  {\boldsymbol{\theta}}) = \frac{\omega_h}{2}\| \mathbf{u}_h(t^k; \boldsymbol{\mu}_i) - \mathbf{\tilde{u}}_h(t^k; \boldsymbol{\mu}_i,  {\boldsymbol{\theta}_{DF}, \boldsymbol{\theta}_D})\|^2 + \frac{1-\omega_h}{2}  \| \tilde{\mathbf{u}}_n(t^k; \boldsymbol{\mu}_i, \boldsymbol{\theta}_E) -  {\mathbf{u}}_n(t^k; \boldsymbol{\mu}_i,  {\boldsymbol{\theta}_{DF}})\|^2
\end{equation}
and $\boldsymbol{\theta} = (\boldsymbol{\theta}_{E}, \boldsymbol{\theta}_{DF}, \boldsymbol{\theta}_{D})$, with $\omega_h \in [0,1]$. The \emph{per-example} loss function (\ref{loss_encoder}) combines the reconstruction error (that is, the error between the FOM solution and the DL-ROM approximation) and the error between the  {intrinsic coordinates} and the output of the encoder. 
This further term allows to enhance the performance of the DL-ROM, as shown in Test 3 of \autoref{sec:4}.

\subsection{Training and Testing Algorithms}

Let us now detail the algorithms through which the training and testing phases of the networks are performed. 

First of all, data normalization and  standardization enhance the training phase of the network by rescaling all the values contained in the dataset to a common frame. For this reason, the inputs and the output of DL-ROM are normalized by applying an affine transformation in order to rescale them in the range $[0, 1]$. In particular, provided a training dataset $X = [\mathbf{x}^1 \, | \, \mathbf{x}^2 \, |  \, \ldots \,  | \,  \mathbf{x}^{N_s}]^T \in \mathbb{R}^{N_s \times M}$, we define
\begin{equation}
\label{max_min}
X_{max} = \max_{i \in \{1, \ldots, N_s\}} \max_{j \in \{1, \ldots, M\}} X_{ij} \quad \textnormal{and} \quad X_{min} = \min_{i \in \{1, \ldots, N_s\}} \min_{j \in \{1, \ldots, M\}} X_{ij}
\end{equation}
so that   data are normalized by applying the following transformation
\begin{equation}
\label{normalization}
X \mapsto \frac{X - X_{min}}{X_{max} - X_{min}}.
\end{equation} \\
Transformation (\ref{normalization}) is applied also to the validation and testing sets, but considering as $X_{max}$ and $X_{min}$ the values computed over the training set. We point out   that the input of the encoder function, the FOM solution $\mathbf{u}_h = \mathbf{u}_h(t^k; \boldsymbol{\mu}_i)$ for a given (time, parameter) instance $(t^k, \boldsymbol{\mu}_i)$, is reshaped in a matrix. In particular, starting from $\mathbf{u}_h \in \mathbb{R}^{N_h}$ we apply the transformation $\mathbf{u}_h^R$=reshape$(\mathbf{u}_h)$ where $\mathbf{u}_h^R \in \mathbb{R}^{N_h^{1/2} \times N_h^{1/2}}$. If $N_h$ is not a square, the input $\mathbf{u}_h$ is zero-padded \citep{goodfellow2016deep}. For the sake of simplicity, we continue to refer to the reshaped FOM solution to as $\mathbf{u}_h$. The inverse reshaping transformation is applied to the output of the last convolutional layer in the decoder function, the ROM approximation. Moreover, we highlight that applying one of the functions (\ref{h_learning})-(\ref{g_learning})-(\ref{encoder}) to the matrix X means applying it row-wise.

The training algorithm referring to the architecture of DL-ROM depicted in \figurename~\ref{architecture_encoder} is reported in Algorithm \ref{training_algorithm}. During the training phase, the optimal parameters of the DL-ROM neural network are found by solving the optimization problem (\ref{minimization_problem_encoder})-(\ref{loss_encoder}) through the back-propagation and ADAM algorithms. 
\begin{algorithm}[t!]
\caption{DL-ROM training algorithm}
\begin{algorithmic}[1]
\Require Parameter matrix $M \in \mathbb{R}^{(n_{\boldsymbol{\mu}} + 1) \times N_s}$, snapshot matrix $S \in \mathbb{R}^{N_h \times N_s}$, training-validation splitting fraction $\alpha$, starting learning rate $\eta$, batch size $N_b$, maximum number of epochs $N_{epochs}$, early stopping criterion, number of minibatches $N_{batches} = (1 - \alpha)N_s/N_b$.
\Ensure  Optimal model parameters $\boldsymbol{\theta}^*= (\boldsymbol{\theta}_E^*, \boldsymbol{\theta}_{DF}^*, \boldsymbol{\theta}_D^*)$.
\vspace{0.3cm}
\State Randomly shuffle $M$ and $S$ \;
\State Split data in $M = [M^{train}, M^{val}]$ and $S = [S^{train}, S^{val}]$ ($M^{val}, S^{val} \in \mathbb{R}^{N_h \times \alpha N_s}$)\;
\State Normalize data in   $M$ and $S$ according to (\ref{normalization})\;
\State Randomly initialize $\boldsymbol{\theta}^0=(\boldsymbol{\theta}_{E}^0, \boldsymbol{\theta}_{DF}^0, \boldsymbol{\theta}_{D}^0)$\;
\State $n_{epochs} = 0$
\While{($\neg$early-stopping \textbf{and} $n_{epochs} \le N_{epochs}$)}
\For{$k = 1 : N_{batches}$}
    \State Sample a minibatch $(M^{batch}, S^{batch}) \subseteq (M^{train}, S^{train})$\;
    \State $S^{batch} =$ reshape$(S^{batch})$\;
    \State $\widetilde{S}^{batch}_n(\boldsymbol{\theta}_{E}^{N_{batches} n_{epochs} + k}) = \mathbf{f}_{n}^E(S^{batch}; \boldsymbol{\theta}_{E}^{N_{batches} n_{epochs} + k})$\;
    \State $S^{batch}_n(\boldsymbol{\theta}_{DF}^{N_{batches} n_{epochs} + k}) = \boldsymbol{\phi}_n^{DF}(M^{batch}; \boldsymbol{\theta}_{DF}^{N_{batches} n_{epochs} + k})$\;
    \State $\widetilde{S}^{batch}_h(\boldsymbol{\theta}_{DF}^{N_{batches} n_{epochs} + k}, \boldsymbol{\theta}_{D}^{N_{batches} n_{epochs} + k}) = \mathbf{f}_{h}^D(S^{batch}_n(\boldsymbol{\theta}_{DF}^{N_{batches} n_{epochs} + k}); \boldsymbol{\theta}_{D}^{N_{batches} n_{epochs} + k})$
    \State $\widetilde{S}^{batch}_h =$ reshape$(\widetilde{S}^{batch}_h)$\;
    \State Accumulate loss (\ref{loss_encoder})  {on $(M^{batch}, S^{batch})$} and compute $\widehat{\nabla}_{\theta} \mathcal{J}$\;
  	\State $\boldsymbol{\theta}^{N_{batches} n_{epochs} + k + 1} = \textnormal{ADAM}(\eta, \widehat{\nabla}_{\theta} \mathcal{J}, \boldsymbol{\theta}^{N_{batches} n_{epochs} + k})$\;
  \EndFor
  \State Repeat instructions 9-13 on $(M^{val}, S^{val})$ with the updated weights $\boldsymbol{\theta}^{N_{batches} n_{epochs} + k + 1}$
  \State Accumulate loss (\ref{loss_encoder}) on $(M^{val}, S^{val})$ to evaluate early-stopping criterion
  \State $n_{epochs} = n_{epochs} + 1$
\EndWhile
\end{algorithmic}
\label{training_algorithm}
\end{algorithm} 

At testing time, the encoder function is instead discarded, that is the DL-ROM architecture is the one shown in \figurename~\ref{architecture_DL-ROM} and the testing algorithm is the one pointed out in Algorithm \ref{testing_algorithm}. The testing phase corresponds to a forward step of the DL-ROM neural network in \figurename~\ref{architecture_DL-ROM}.
\begin{algorithm}[ht!]
\caption{DL-ROM testing algorithm}
\begin{algorithmic}[1]
\Require Testing parameter matrix $M^{test} \in \mathbb{R}^{(n_{\boldsymbol{\mu}} + 1) \times (N_{test} N_t)}$ and optimal model parameters $(\boldsymbol{\theta}_{DF}^*, \boldsymbol{\theta}_D^*)$.
\Ensure ROM approximation matrix $\widetilde{S}_h \in \mathbb{R}^{N_h \times (N_{test} N_t)}$.
\vspace{0.3cm}
\State Load $\boldsymbol{\theta}_{DF}^*$ and $\boldsymbol{\theta}_D^*$\;
\State $S_n(\boldsymbol{\theta}_{DF}^*) = \boldsymbol{\phi}_n^{DF}(M^{test}; \boldsymbol{\theta}_{DF}^*)$\;
\State $\widetilde{S}_h(\boldsymbol{\theta}_{DF}^*, \boldsymbol{\theta}_{D}^*) = \mathbf{f}_{h}^D(S_n(\boldsymbol{\theta}_{DF}^*); \boldsymbol{\theta}_{D}^*)$
\State $\widetilde{S}_h=$ reshape$(\widetilde{S}_h)$\;
\end{algorithmic}
\label{testing_algorithm}
\end{algorithm}

\section{Numerical results}
\label{sec:4}

In this section, we report the numerical results obtained by applying  the proposed DL-ROM technique to three parametrized, time-dependent PDE problems, namely {\em (i)} Burgers equation, {\em (ii)} a linear transport equation, and {\em (iii)} a coupled  PDE-ODE system arising from cardiac electrophysiology, namely the monodomain equation; this latter is a system of time dependent, nonlinear equations, whose solutions feature a traveling wave behavior. For the time being, we deal with problems set in $d=1$ (spatial) dimension featuring up to $n_{\mu} = 2$ parameters; we will consider the extension to differential problems in $d=2$ and $d=3$ in a forthcoming publication. For this reason, our focus is now on the numerical accuracy of our DL-ROM technique rather than on its computational efficiency and, therefore, on its comparison with linear ROMs such as the RB method featuring linear (possibly, piecewise linear) trial manifolds built through POD. 

To evaluate the performance of DL-ROM we rely on the loss function (\ref{loss_encoder}) and on the following error indicator
\begin{equation}
\epsilon_{rel} = \frac{1}{N_{test}} \sum_{i  = 1}^{N_{test}} \left(\displaystyle \frac{\sqrt{ \sum_{k=1}^{N_t} || \mathbf{u}^k_h(\boldsymbol{\mu}_{test,i}) - \mathbf{\tilde{u}}^k_h(\boldsymbol{\mu}_{test,i}) ||^2}}{\sqrt{\sum_{k=1}^{N_t} || \mathbf{u}_h^k(\boldsymbol{\mu}_{test,i}) ||^2}} \right).
\label{relative_error}
\end{equation}
We implement the neural network required by our DL-ROM technique by means of the Tensorflow deep learning framework \citep{abadi2016tensorflow} and the numerical simulations are performed on a workstation equipped with an Nvidia GeForce GTX 1070 8 GB GPU.

\subsection{Test 1: Burgers Equation}
Let us consider the parametrized one-dimensional nonlinear Burgers equation
\begin{equation}
\label{Burgers}
\begin{cases}
\displaystyle
\frac{\partial u}{\partial t} + u \frac{\partial u}{\partial x} - \frac{1}{\mu} \frac{\partial^2 u}{\partial x^2} = 0, \quad & (x, t) \in (0, L) \times (0, T) \\
u(0, t) = 0, \quad & t \in (0, T) \\
u(L, t) = 0, \quad & t \in (0, T) \\
u(x, 0) = u_0(x), \quad & x \in (0, L),
\end{cases}
\end{equation} \\
where
\begin{equation*}
u_0(x) = \frac{x}{1 + \sqrt{1 / A_0} \exp( \mu x^2 / 4)},
\end{equation*}
with $A_0 = \exp(\mu / 8)$, $L = 1$ and $T =2$. System (\ref{Burgers}) has been discretized in space by means of linear finite elements, with $N_h =256$ grid points, and in time by means of the Backward Euler scheme, with $N_t = 100$ time instances. The parameter space, to which belongs the single ($n_{\mu} = 1$) parameter, is given by $\mathcal{P} = [100, 1000]$. We consider $N_{train} = 20$ training-parameter instances uniformly distributed over $\mathcal{P}$ and $N_{test} = 19$ testing-parameter instances, each of them corresponding to the midpoint between two consecutive training-parameter instances. 

The configuration of the DL-ROM neural network used for this test case is the following. We choose a 12-layers DFNN equipped with 50 neurons per hidden layer and $n$ neurons in the output layer, where $n$ corresponds to the dimension of the reduced  {trial} manifold. The architectures of the encoder and decoder functions are instead  reported in \tablename s \ref{table_transposed_convolutional_layers_encoder} and \ref{table_transposed_convolutional_layers}, and are similar to the ones used in \citep{carlberg2018model}.
\begin{table}[ht]
{\small
\begin{center}
\begin{tabular}{|c|c|c|c|c|c|c|}
\hline
Layer & Input Dimension & Output Dimension & Kernel Size & $\#$ of Filters & Stride & Padding \\
\hline
1 & & & [5, 5] & 8 & 1 & SAME \\
\hline
2 & & & [5, 5] & 16 & 2 & SAME \\
\hline
3 & & & [5, 5] & 32 & 2 & SAME \\
\hline
4 & & &[5, 5] & 64 & 2 & SAME \\
\hline
5 & $N_h$ & 256 & & & & \\
\hline
6 & 256 & $n$ & & & &\\
\hline
\end{tabular}
\end{center}
}
\vspace{-0.25cm}
\caption{\textit{Test 1}: Attributes of convolutional layers and dense layers in the encoder $\mathbf{f}_n^E$.}
\label{table_transposed_convolutional_layers_encoder}
\end{table}
\begin{table}[ht]
{\small
\begin{center}
\begin{tabular}{|c|c|c|c|c|c|c|}
\hline
Layer & Input dimension & Output dimension & Kernel size & $\#$ of filters & Stride & Padding \\
\hline
1 & $n$ & 256 & & & &\\
\hline
2 & 256 & $N_h$ & & & &\\
\hline
3 & & & [5, 5] & 64 & 2 & SAME \\
\hline
4 & & & [5, 5] & 32 & 2 & SAME \\
\hline
5 & & & [5, 5] & 16 & 2 & SAME \\
\hline
6 & & & [5, 5] & 1 & 1 & SAME \\
\hline
\end{tabular}
\end{center}
}
\vspace{-0.25cm}
\caption{\textit{Test 1}: Attributes of dense layers and transposed convolutional layers in the decoder $\mathbf{f}_h^D$.}
\label{table_transposed_convolutional_layers}
\end{table} 

\begin{figure}[ht!]
\centering
\includegraphics[scale=0.45]{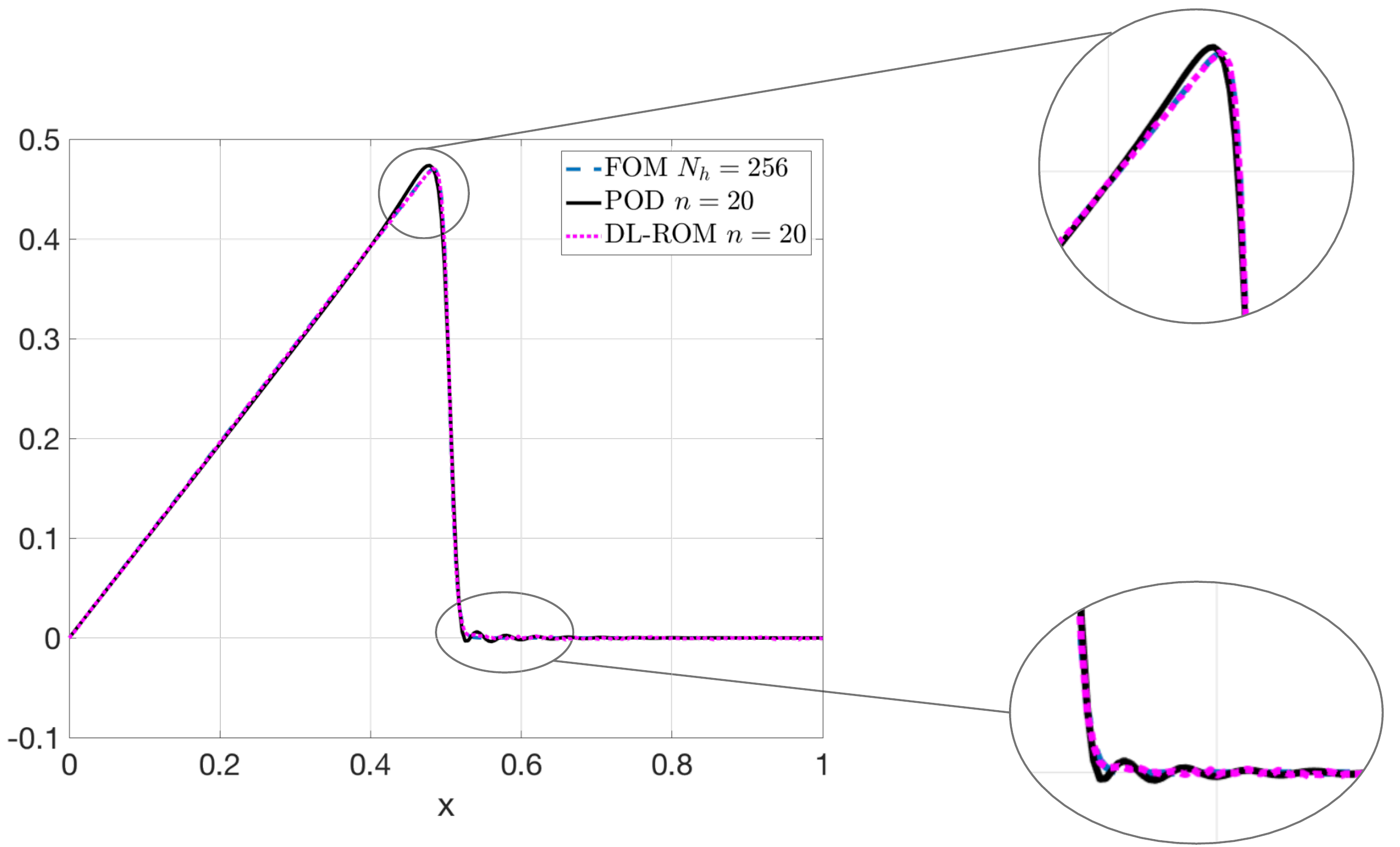}
\vspace{-0.15cm}
\caption{\textit{Test 1}: FOM, optimal-POD  and DL-ROM solutions for the testing-parameter instance $\mu_{test} = 976.32$ at $t = 0.02$, with $n=20$.}
\label{comparison_burgers}
\end{figure}

Problem (\ref{Burgers}) does not represent a remarkably challenging task for linear ROM, indeed by considering for example POD and by applying it to the snapshot matrix (the latter built by collecting the solution of (\ref{Burgers}) for $N_s = N_{train}N_t$ training-parameter instances) it is sufficient to assemble a linear trial manifold of dimension 20 in order to capture more than the 99.99$\%$ of the energy of the system \citep{san2018neural, quarteroni2016reduced}. In order to assess the performance of our DL-ROM technique, we compute the DL-ROM solution by fixing the dimension of the nonlinear trial manifold to $n=20$. In \figurename~\ref{comparison_burgers} we show the DL-ROM and the optimal-POD reconstructions, along with the FOM solution, for the time instance $t = 0.02$ and for the testing-parameter instance $\mu_{test} = 976.32$, the testing value of $\mu$ for which the reconstruction task results to be the most difficult both for POD and DL-ROM, being the diffusion term in (\ref{Burgers}) smaller and the solution closer to the one of a purely hyperbolic system. In particular, for $\mu_{test} = 976.32$, employing the DL-based ROM technique presented in this work allows us to halve the error indicator $\epsilon_{rel}$ associated to the optimal-POD approximation of the FOM solution. Referring to \figurename~\ref{comparison_burgers}, the DL-ROM reconstruction is more accurate than the optimal-POD one, indeed it mostly fits the FOM solution, even in correspondence of its maximum, as shown in the zooms of \figurename~\ref{comparison_burgers}. Moreover, it does not introduce oscillations where a large gradient of the FOM solution is observed, as it happens instead by employing POD.

\begin{figure}[ht!]
\centering
\includegraphics[scale=0.145]{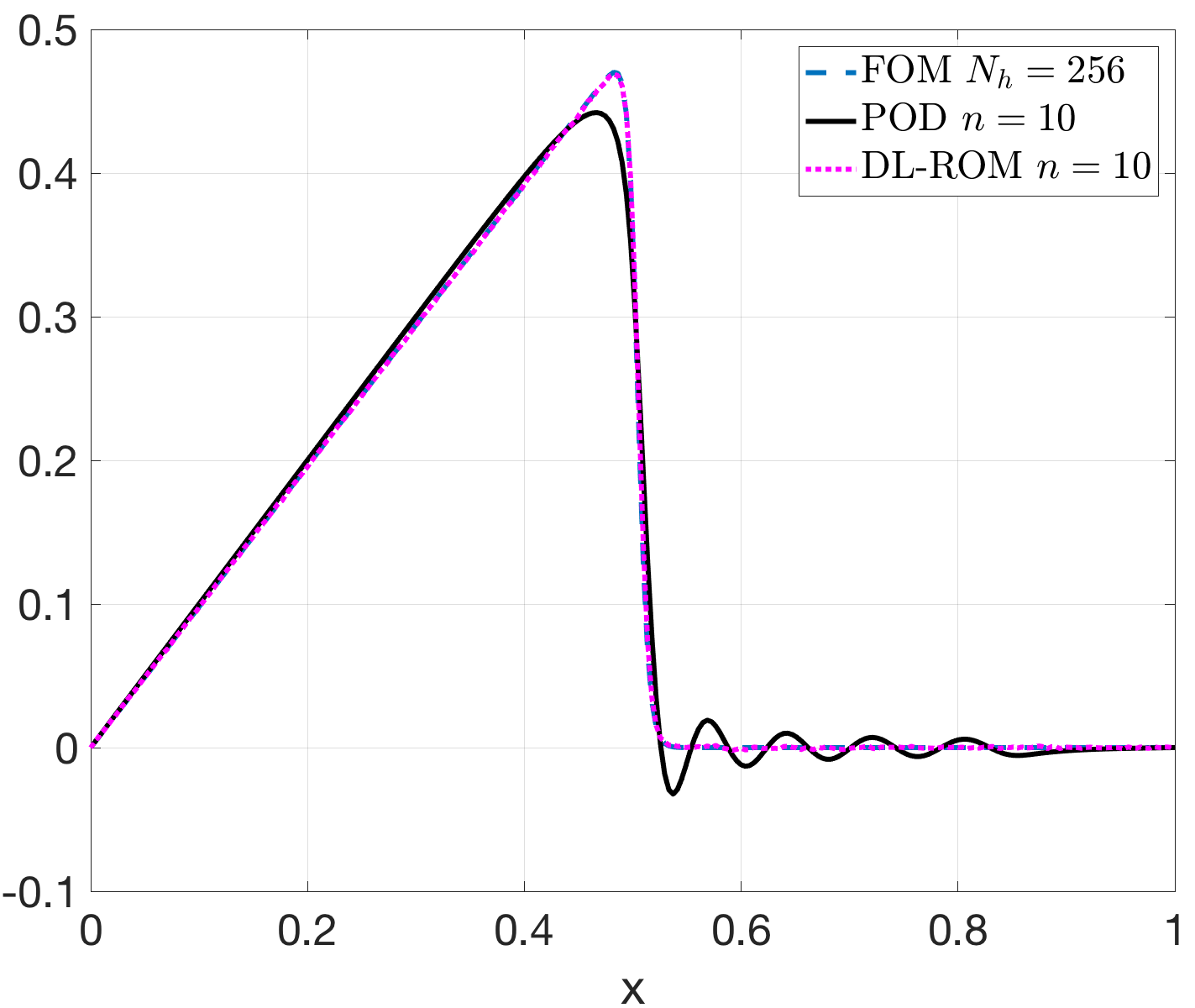}
\vspace{-0.15cm}
\caption{\textit{Test 1}: FOM,  optimal-POD  and DL-ROM solutions for the testing-parameter instance $\mu_{test} = 976.32$ at $t = 0.02$, with $n=10$.}
\label{comparison_burgers_p=10}
\end{figure}

In \figurename~\ref{comparison_burgers_p=10} we show the same comparison of \figurename~\ref{comparison_burgers} but this time considering both for POD and DL-ROM a reduced dimension $n=10$. The difference in terms of accuracy provided by the two approaches is even more striking in this case.

\begin{figure}[ht!]
\centering
\includegraphics[scale=0.135]{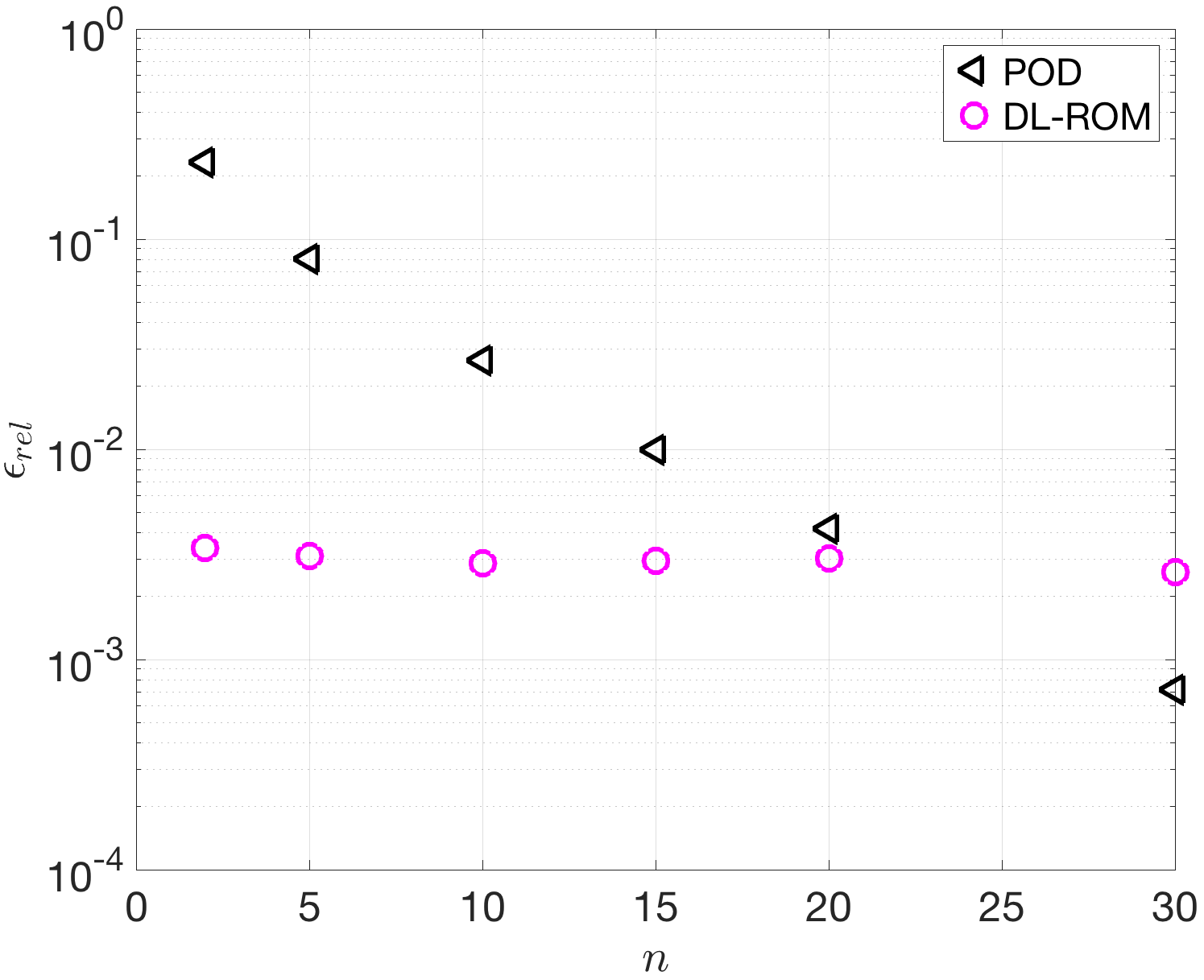}
\vspace{-0.15cm}
\caption{\textit{Test 1}: Error indicator $\epsilon_{rel}$ vs. $n$ on the testing set.}
\label{1D_burgers_convergence}
\end{figure}

Finally, in \figurename~\ref{1D_burgers_convergence} we highlight the accuracy properties of both the DL-ROM and POD techniques by displaying the behavior of the error indicator $\epsilon_{rel}$, defined in (\ref{relative_error}), with respect to the dimension $n$ of the corresponding reduced {trial} manifold. For $n < 20$ the DL-ROM approximation is more accurate than the one provided by POD, and only for $n = 20$ the two techniques provide almost the same accuracy.

\subsection{Test 2: Linear Transport Equation}

We consider two tests for this set of parametrized differential models.

\subsubsection*{Test 2.1: $n_{\mu} = 1$}

First, we consider the parametrized one-dimensional linear transport equation
\begin{equation}
\label{Transport}
\begin{cases}
\displaystyle
\frac{\partial u}{\partial t} + \mu \frac{\partial u}{\partial x} = 0, \quad & (x, t) \in \mathbb{R} \times (0, T) \\
u(x, 0) = u_0(x), \quad & x \in \mathbb{R},
\end{cases}
\end{equation} \\
whose exact solution is $u(x, t) = u_0(x - \mu t)$. We set $u_0(x) = (1/\sqrt{2 \pi \sigma}) e^{-x^2/2 \sigma}$ and $T = 1$. 

The parameter (here $n_{\mu} = 1$) represents the velocity of the travelling wave and the parameter space is given by $\mathcal{P} = [0.775, 1.25]$. The dataset is built by uniformly sampling the exact solution in the domain $(0, L) \times (0, T)$, with $L = 1$, and by considering $N_h = 256$ degrees of freedom in the space discretization and $N_t = 200$ time instances in the time one. We consider $N_{train} = 20$ training-parameter instances uniformly distributed in the parameter space $\mathcal{P}$ and $N_{test} = 19$ testing-parameter instances such that $\mu_{test,i}= (	\mu_{train,i}+ \mu_{train,i+1})/2$, for $i = 1, \ldots, N_{test}$. This test case, and more in general hyperbolic problems, are examples in which the use of a linear approach to ROM generally yields poor performance in terms of accuracy. Indeed, the dimension of the linear trial manifold must be very large, if compared to the dimension of the solution manifold, in order to capture the variability of the FOM solution over the parameter space $\mathcal{P}$. We set $\sigma = 10^{-4}$ in order to assess the performance of DL-ROM in a scenario which is still remarkably challenging for ROM on linear trial manifolds.

\begin{figure}[ht!]
\centering
\includegraphics[scale=0.325]{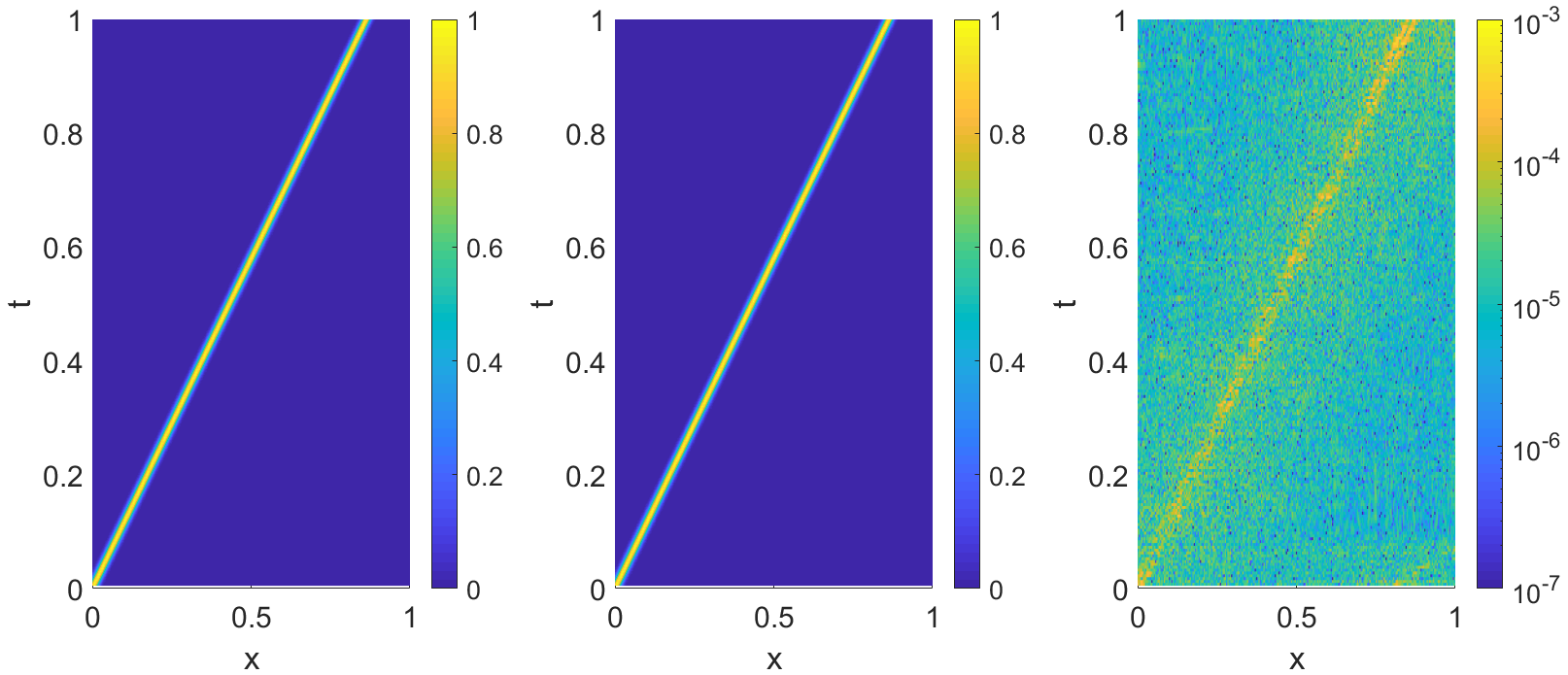}
\vspace{-0.15cm}
\caption{\textit{Test 2.1}: Exact solution (left), DL-ROM solution with $n = 2$ (center) and relative error $\boldsymbol{\epsilon}_k$ (right) for the testing-parameter instance $\mu_{test} = 0.8625$ in the space-time domain.}
\label{comparison_transport}
\end{figure}

\figurename~\ref{comparison_transport} shows the exact solution, which here plays the role of the FOM solution, and the DL-ROM one for the testing-parameter instance $\mu_{test} = 0.8625$; here, we set the dimension of the nonlinear trial manifold to $n = 2$, equal to the dimension of the solution manifold $n_{\mu} + 1$. Moreover, in \figurename~\ref{comparison_transport} we highlight the relative error $\boldsymbol{\epsilon}_k \in \mathbb{R}^{N_h}$, for $k = 1, \ldots, N_t$, associated to a given $\boldsymbol{\mu}_{test} \in \mathcal{P} \subset \mathbb{R}^{n_{\mu}}$ (in this case $n_{\mu} = 1$), 
defined as
\begin{equation}
\boldsymbol{\epsilon}_k = \displaystyle \frac{ | \mathbf{u}^k_h(\boldsymbol{\mu}_{test}) - \mathbf{\tilde{u}}^k_h(\boldsymbol{\mu}_{test}) |}{\sqrt{\sum_{k=1}^{N_t} || \mathbf{u}^k_h(\boldsymbol{\mu}_{test}) ||^2}},
\label{error}
\end{equation}
which widens in proximity of the spike of the exact solution.

\begin{figure}[ht]
\centering
\includegraphics[scale=0.15]{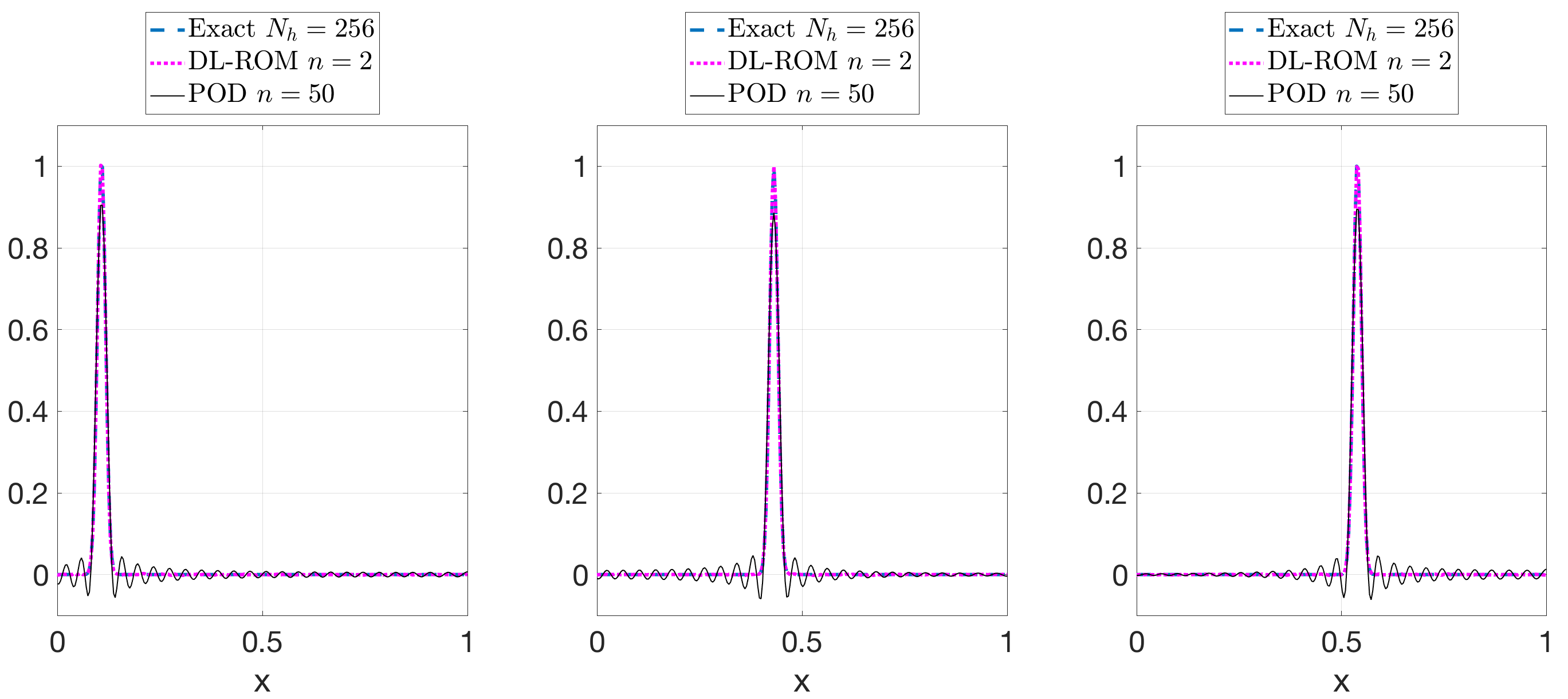}
\caption{\textit{Test 2.1}: Exact, DL-ROM and POD solutions for the testing-parameter instance $\mu_{test} = 0.8625$ at $t = 0.125, 0.5$ and 0.625.}
\label{comparison_transport_time}
\end{figure}

In \figurename~\ref{comparison_transport_time} we report the exact solution and the DL-ROM one, obtained by setting $n = 2$, for three particular time instances. In order to compare the performance of the proposed nonlinear ROM with a linear approach, we perform the POD on the snapshot matrix and show, for the same testing-parameter instance, the optimal POD-reconstruction, i.e. the projection of the FOM (exact) solution onto the POD basis, in \figurename~\ref{comparison_transport_time}. For example, by considering $n = 2$, the error indicator, defined in (\ref{relative_error}), is $\epsilon_{rel} = 8.74 \cdot 10^{-3}$. By considering a linear ROM technique instead, even by considering a reduced trial manifold of dimension $n = 50$, built by means of the POD, the reconstructed solution presents spurious oscillations which result in a poor approximation of the FOM solution (see \figurename~\ref{comparison_transport_time}). Indeed, in order to achieve the same accuracy obtained through DL-ROM over the testing set one has to select 90 basis functions, i.e. a linear trial manifold of dimension $n = 90$.

\begin{figure}[ht]
\centering
\includegraphics[scale=0.14]{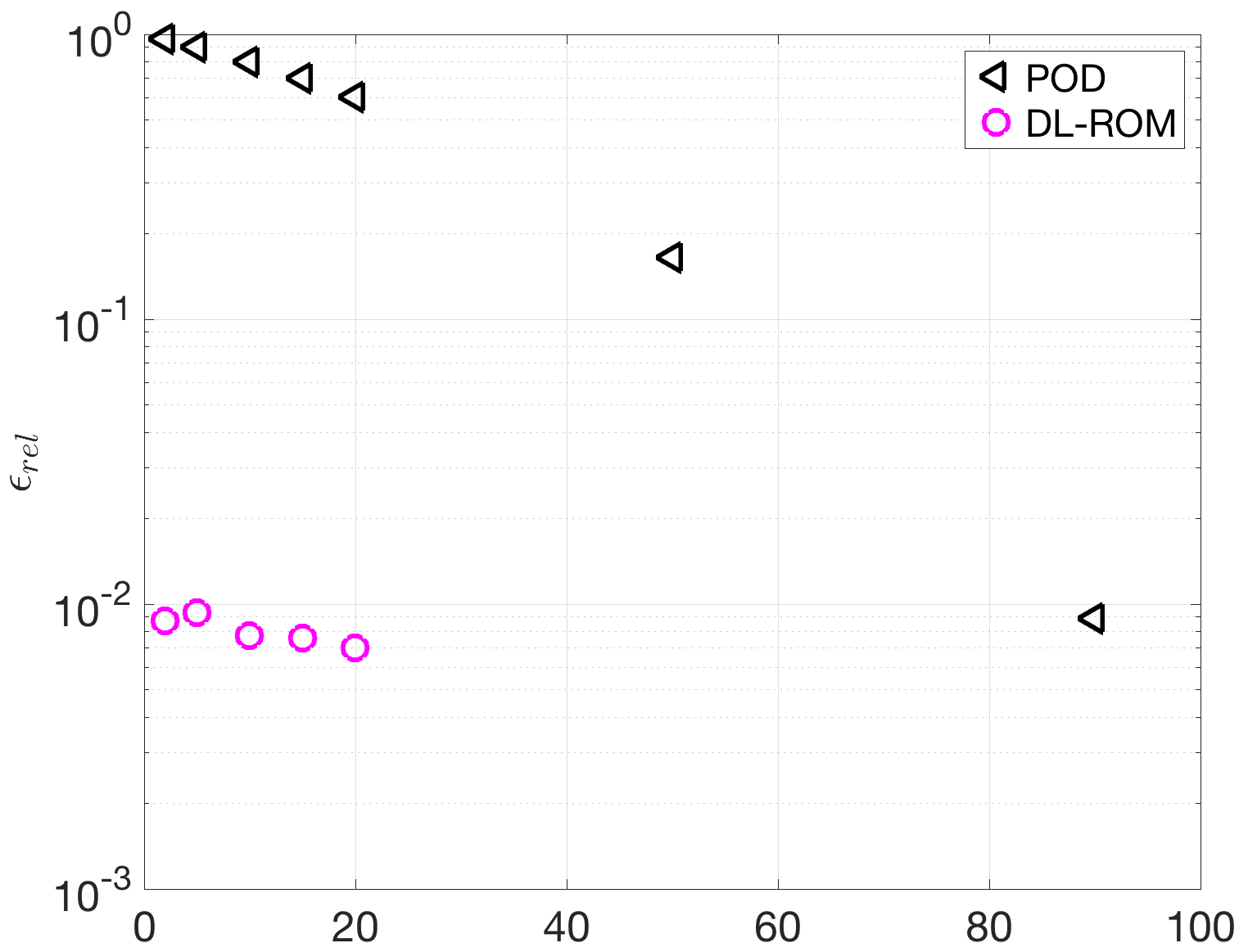}
\caption{\textit{Test 2.1}: Error indicator $\epsilon_{rel}$ vs. $n$ on the testing set.}
\label{relative_error_transport}
\end{figure}

\figurename~\ref{relative_error_transport} shows the behavior of the error indicator (\ref{relative_error}) with respect to the reduced dimension $n$. By increasing the dimension of the nonlinear trial manifold there is a slight improvement of the performance of the DL-ROM neural network, i.e. the error indicator decreases. This improvement is not particularly {relevant} because by increasing $n$, the number of parameters of the DL-ROM neural network, i.e. weights and biases, is increased by a limited {quantity}. In this way the approximation capability of the neural network remains almost the same and so  {does} the error indicator (\ref{relative_error}).

\begin{rem} \textit{(Hyperparameters Tuning)}. The hyperparameters of the DL-ROM neural network are tuned by evaluating the loss function over the validation set and by setting each of them equal to the value minimizing the generalization error on the validation set. In particular, we show the tests performed to choose the size of the (transposed) convolutional kernels in the (decoder) encoder function, the number of hidden layers in the feedforward neural network and the number of neurons for each hidden layer. The hyperparameters evaluation starts from the default configuration in \tablename~\ref{1D_transport_starting_configuration}.

\begin{table}[ht]
\begin{center}
\begin{tabular}{|c|c|c|}
\hline
Kernel Size & $\#$ Hidden Layers & $\#$ Neurons \\
\hline
[3, 3] & 1 & 50 \\
\hline
\end{tabular}
\end{center}
\vspace{-0.25cm}
\caption{\textit{Test 2.1}: Starting configuration of DL-ROM.}
\label{1D_transport_starting_configuration}
\end{table}

Then, the best values are found  {iteratively} by studying the impact of the variation of a single hyperparameter at a time on the validation loss.  {Once the best value of a hyperparameter  is found, this value replaces the default value from that point on}. For each hyperparameter the tuning is performed in a range of values for which the training of the  network is affordable regarding computational costs. 
\begin{figure}[ht]
\centerline{
\includegraphics[width=0.33\textwidth]{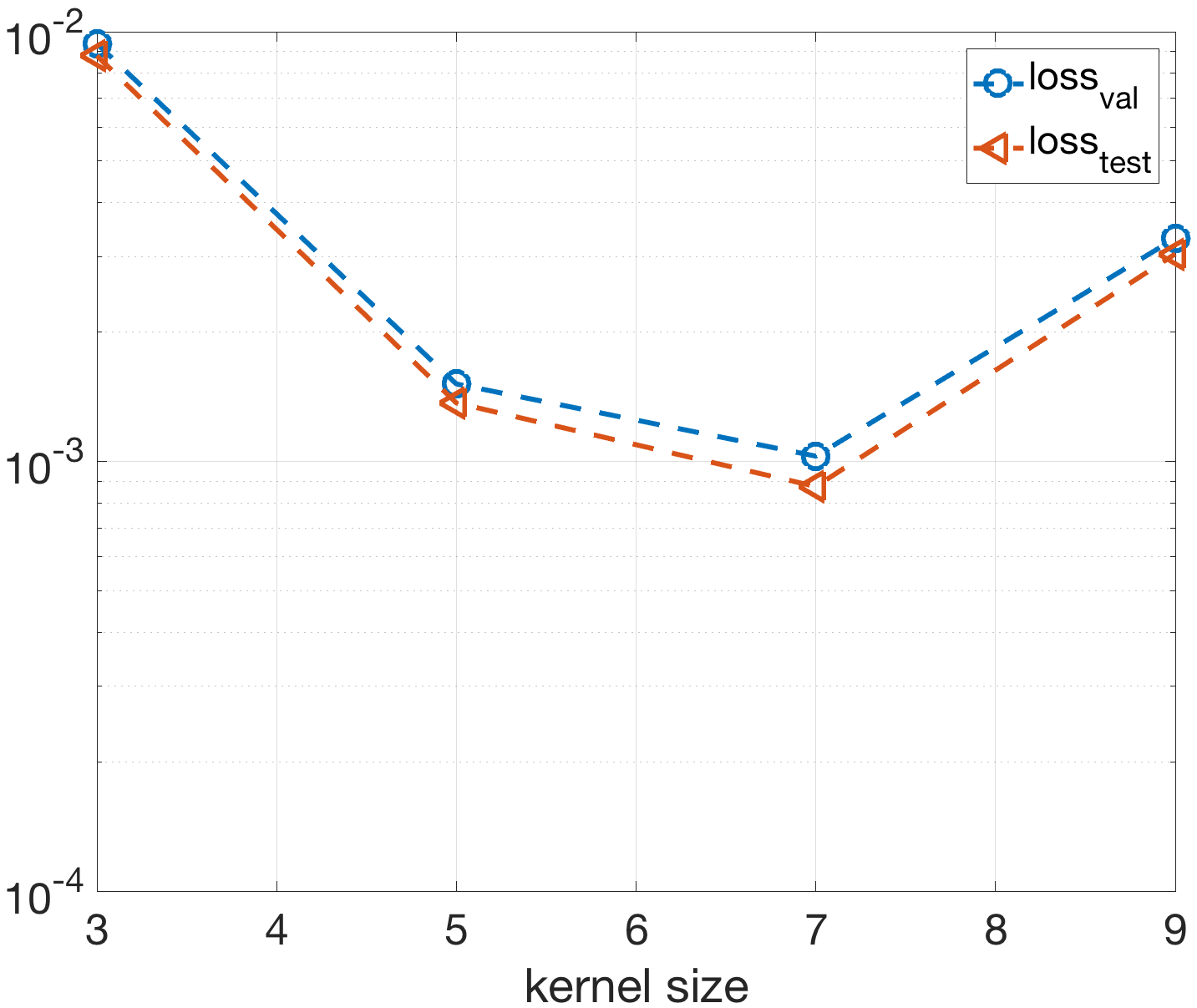}
\includegraphics[width=0.33\textwidth]{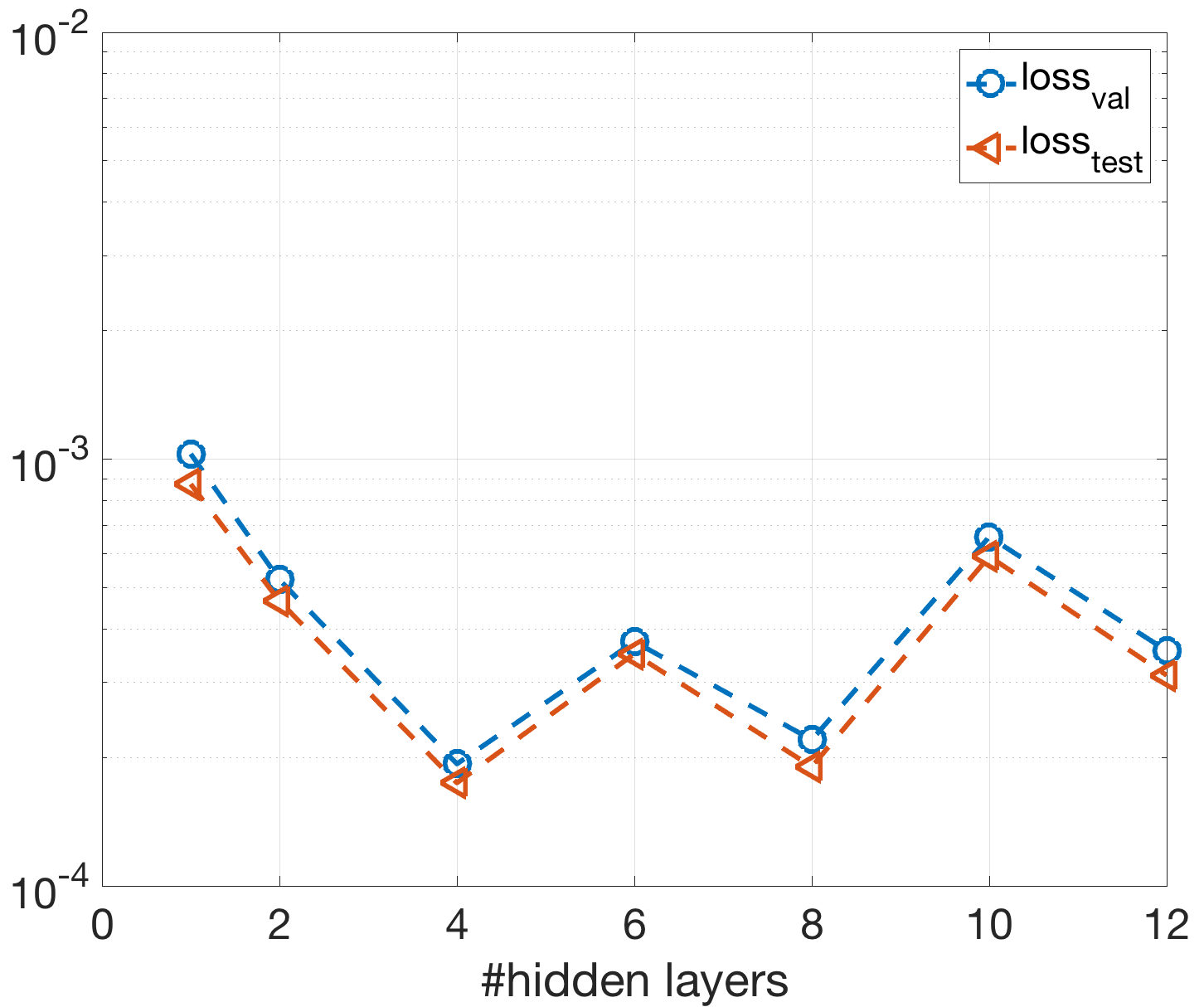}
\includegraphics[width=0.33\textwidth]{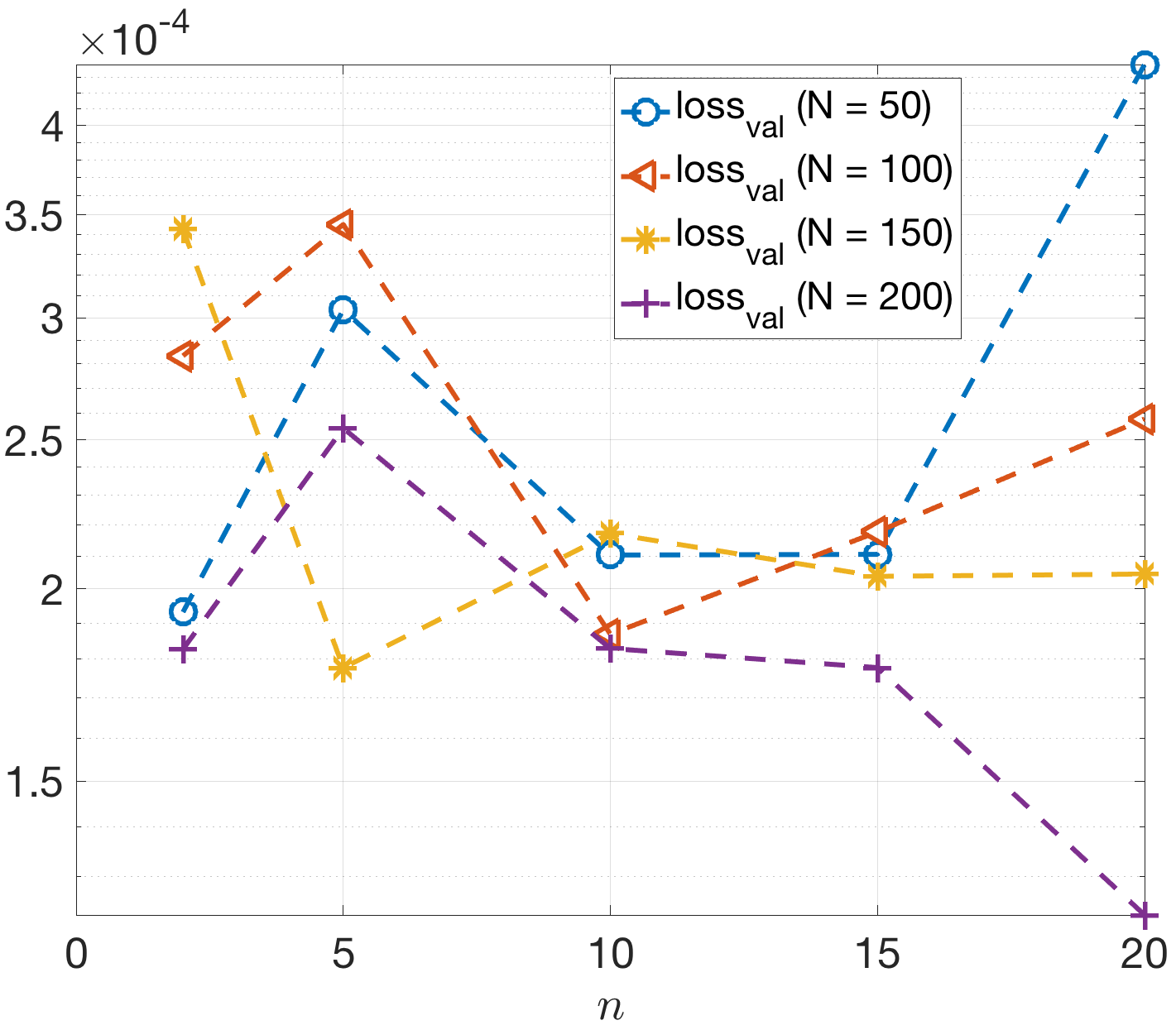}
}
\caption{\textit{Test 2.1}: Impact of the kernel size (left), the number of hidden layers (center) and the number of neurons (right) on the validation and testing loss.}
\label{kernel_size}
\end{figure}
In \figurename~\ref{kernel_size}, we show the impact of the size of the convolutional kernels on the loss over the validation and testing sets, the number of hidden layers in the feedforward forward neural network and the number of neurons in each hidden layer by varying the reduced dimension in order to find the best value of such hyperparameter over $n$. The final configuration of the DL-ROM neural network is the one provided in \tablename~\ref{1D_transport_final_configuration}.
\begin{table}[ht]
\begin{center}
\begin{tabular}{|c|c|c|}
\hline
Kernel Size & $\#$ Hidden Layers & $\#$ Neurons \\
\hline
[7, 7] & 4 & 200 \\
\hline
\end{tabular}
\end{center}
\vspace{-0.2cm}
\caption{\textit{Test 2.1}: Final configuration of DL-ROM.}
\label{1D_transport_final_configuration}
\end{table}
\end{rem}

\subsubsection*{Test 2.2: $n_{\mu} = 2$}

Here we consider again the parametrized one-dimensional transport equation 
\begin{equation}
\label{Transport_2params}
\begin{cases}
\displaystyle
\frac{\partial u}{\partial t} + \frac{\partial u}{\partial x} = 0, \quad & (x, t) \in \mathbb{R} \times (0, T) \\
u(x, 0) = u_0(x), \quad & x \in \mathbb{R}.
\end{cases}
\end{equation} 
The exact solution of (\ref{Transport_2params}) is $u(x, t) = u_0(x - t; \boldsymbol{\mu})$ but this time we set the initial datum equal to
\begin{equation}
u_0(x ; \boldsymbol{\mu}) =
\begin{cases} 
0, \quad &  \textnormal{if} \; x < \mu_1 \\
\mu_2, \quad &  \textnormal{if} \; x \geq \mu_1, \\
\end{cases}
\label{initial_data}
\end{equation}
where $\boldsymbol{\mu}=[\mu_1, \mu_2]^T$. The $n_{\mu}=2$ parameters belong to the parameter space $\mathcal{P} = \mathcal{P}_{\mu_1} \times \mathcal{P}_{\mu_2} = [0.025, 0.25] \times [0.5, 1]$. We build the dataset by uniformly sampling the exact solution in the domain $(0, L) \times (0, T)$, with $L = 1$ and $T=1$, and by considering $N_h = 256$ grid points for the space discretization and $N_t = 100$ time instances for the time one. We collect, both for $\mu_1$ and $\mu_2$, $N_{train} = 21$ training-parameter instances uniformly distributed in the parameter space $\mathcal{P}$ and $N_{test} = 20$ testing-parameter instances, selected as in the other test cases. Equation (\ref{Transport_2params}), completed  {with} the initial datum (\ref{initial_data}), stands as one of the most challenging problems for linear ROM techniques because of the difficulty to accurately reconstruct the jump discontinuity of the exact solution as a linear combination of basis functions computed from the snapshots, for a testing-parameter instance. The architecture of the DL-ROM neural network used here is the one presented in the Test 2.1.

\begin{figure}[ht!]
\centering
\includegraphics[scale=0.35]{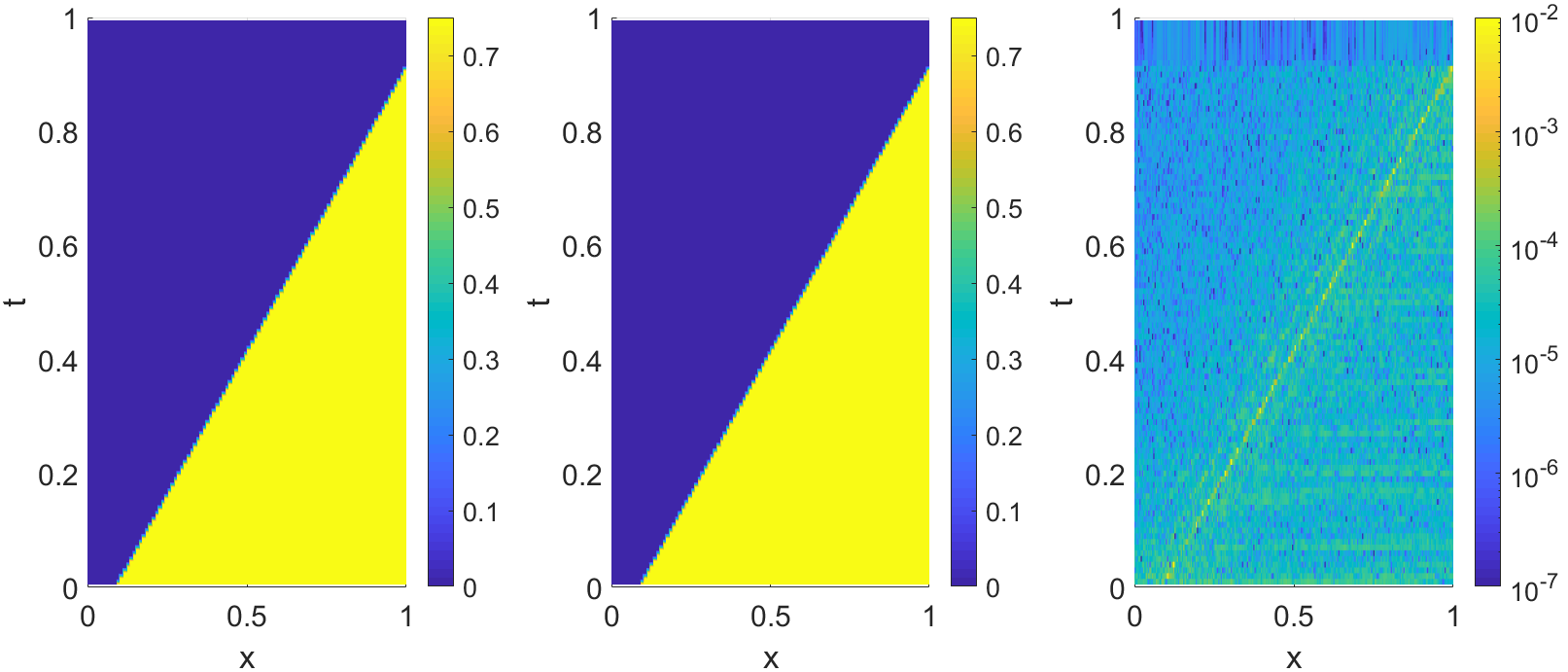}
\caption{\textit{Test 2.2}: Exact solution (left), DL-ROM solution with $n = 3$ (center) and relative error $\boldsymbol{\epsilon}_k$ (right) for the testing-parameter instance $\boldsymbol{\mu}_{test} = (0.154375, 0.6375)$ in the space-time domain.}
\label{comparison_1D_transport_2params}
\end{figure}

In \figurename~\ref{comparison_1D_transport_2params} we show the exact solution, which here again plays the role of the FOM solution, and the DL-ROM one, obtained by setting $n=3$, equal to the dimension of the solution manifold $n_{\mu} + 1$, for the testing-parameter instance $\boldsymbol{\mu}_{test} = (0.154375, 0.6375)$, along with the {relative error} $\boldsymbol{\epsilon}_k$, defined in (\ref{error}), which is larger near the jump of the FOM solution.

\begin{figure}[ht!]
\centering
\includegraphics[scale=0.15]{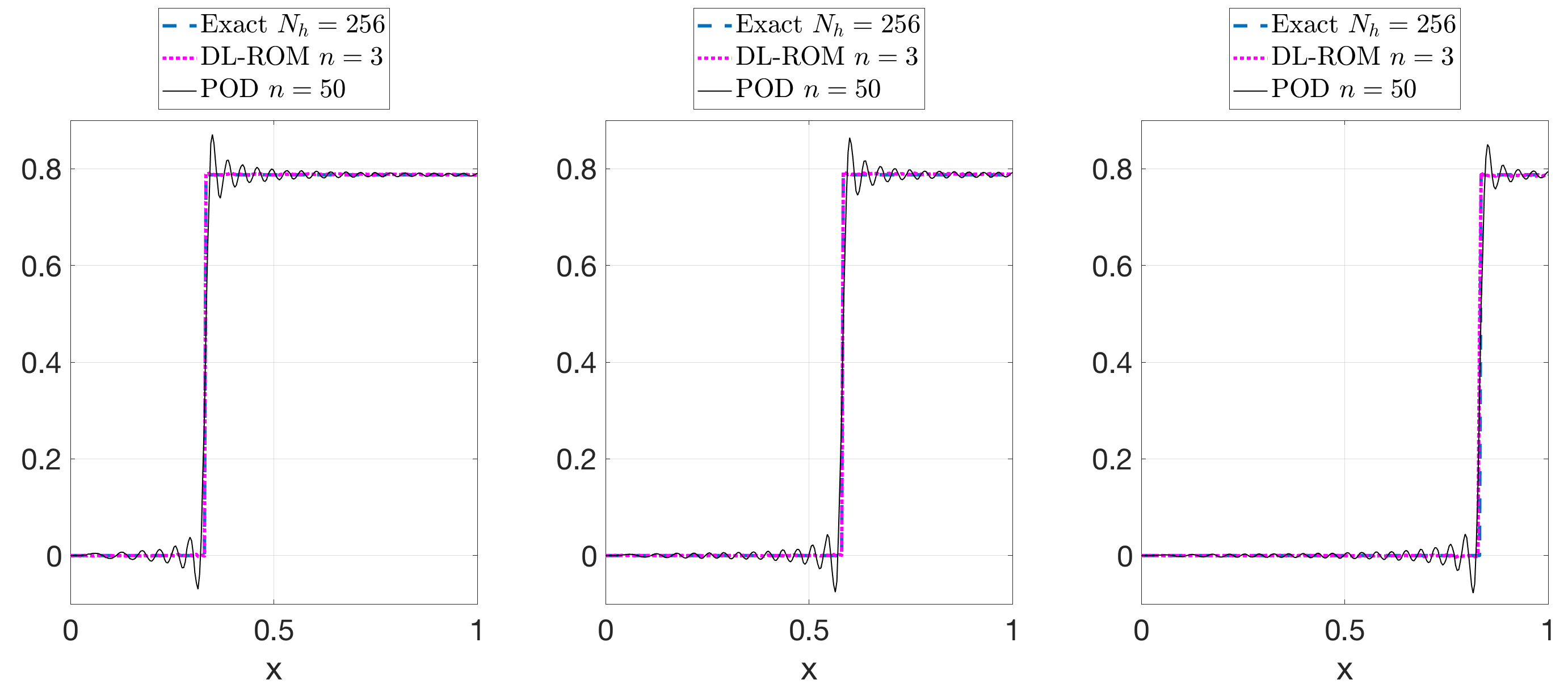}
\caption{\textit{Test 2.2}: Exact, DL-ROM and optimal-POD solutions for the testing-parameter instance $\boldsymbol{\mu}_{test} = (0.154375, 0.6375)$ at $t = 0.245, 0.495$ and 0.745.}
\label{comparison_1D_transport_time_2params}
\end{figure}

In \figurename~\ref{comparison_1D_transport_time_2params} we report the DL-ROM and optimal-POD reconstructions, together with the FOM solution, for the time instances $t = 0.245, 0.495$ and 0.745, and the testing-parameter instance $\boldsymbol{\mu}_{test} = (0.154375, 0.6375)$. The dimension of the reduced manifolds are $n = 3$ and $n = 50$ for the DL-ROM and POD techniques, respectively. By considering a linear ROM technique, even by setting the dimension of the reduced manifold equal to $n = 50$, the reconstructed solution presents spurious oscillations which lead to a poor approximation of the FOM solution. Moreover, the optimal-POD solution is not able to   fit the discontinuity of the FOM solution in a sharp way. These oscillations are significantly mitigated by the use of our DL-ROM and the jump discontinuity is accurately fit by the DL-ROM solution, as shown in \figurename~\ref{comparison_1D_transport_time_2params}.

\begin{figure}[ht!]
\centering
\includegraphics[scale=0.14]{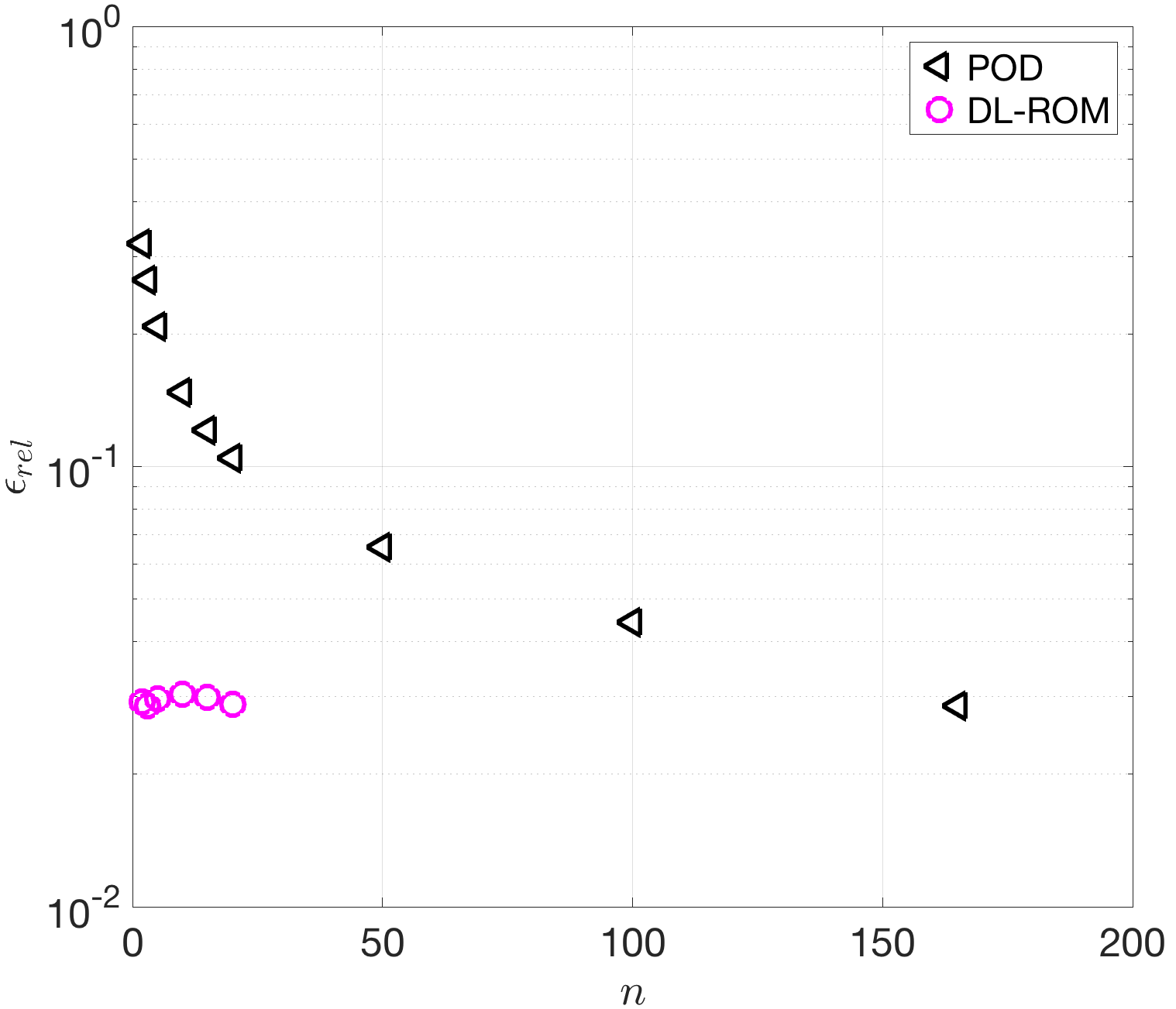}
\caption{\textit{Test 2.2}: Error indicator $\epsilon_{rel}$ vs. $n$ on the testing set.}
\label{relative_error_transport_2params}
\end{figure}

Finally, in \figurename~\ref{relative_error_transport_2params} we highlight the accuracy properties of both the DL-ROM and POD techniques. In particular, the same conclusions observed in Test 2.1, namely those  {regarding} the behaviour of the error indicator (\ref{relative_error}) with respect to the reduced dimension $n$, still hold. The developed DL-ROM technique allows us to obtain a value for the error indicator equal to $\epsilon_{rel} = 2.85 \cdot 10^{-2}$ with $n = 3$, which instead is achieved by POD only by selecting 165 basis functions, i.e. by building a linear trial manifold of dimension $n = 165$.

\subsection{Test 3: Monodomain Equation}

\begin{figure}[ht!]
\centering
\includegraphics[scale=0.145]{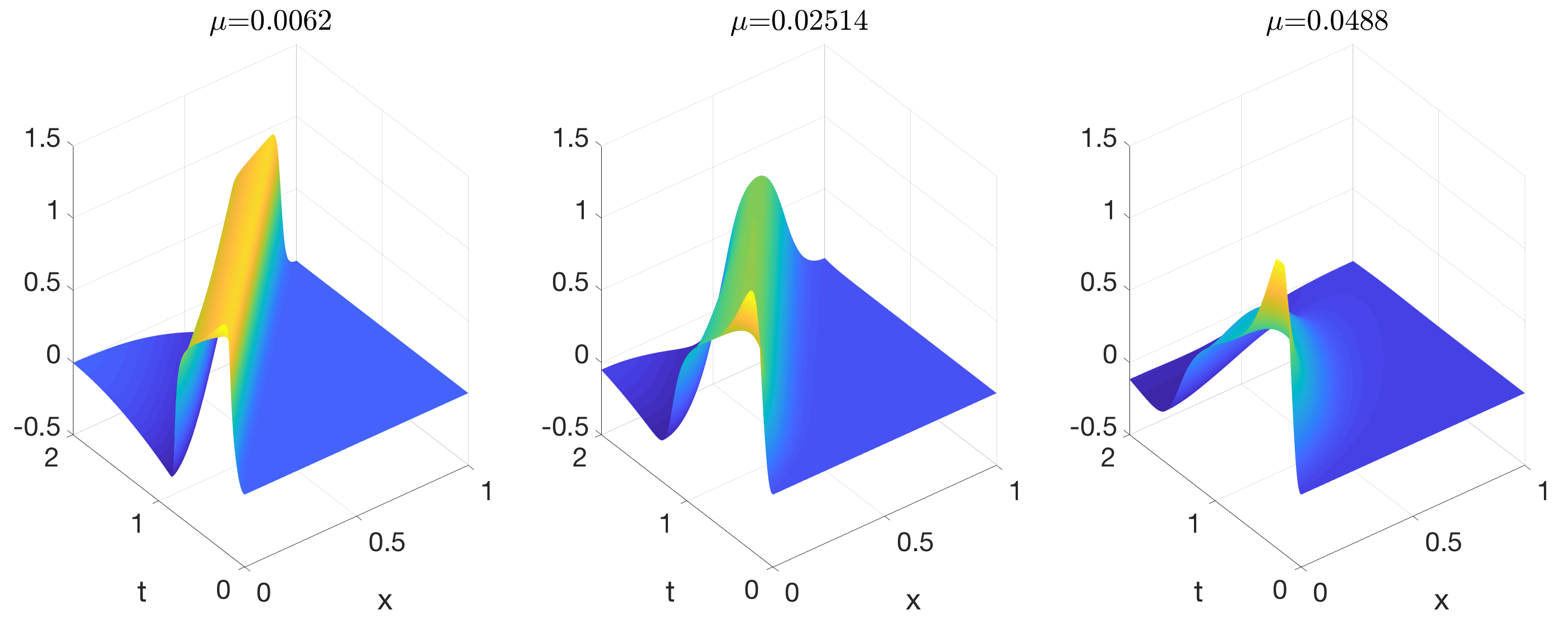}
\caption{\textit{Test 3}: FOM solutions for different testing-parameter instances.}
\label{fom_solution_1Dmonodomain}
\end{figure}

We now consider the following one-dimensional coupled PDE-ODE nonlinear system
\begin{equation}
\label{1DMonodomain}
\begin{cases}
\vspace{0.1cm}
\displaystyle \mu \frac{\partial u}{\partial t} - \mu^2 \frac{\partial^2 u}{\partial x^2} + u(u-0.1)(u-1) + w = 0, \quad & (x, t) \in (0,L) \times (0,T) \\
\vspace{0.1cm}
\displaystyle \frac{d w}{d t} + (\gamma w - \beta u)=0, \quad & (x, t) \in (0,L) \times (0,T) \\
\vspace{0.1cm}
\displaystyle \frac{\partial u}{\partial x}(0,t) = 50000 t^3 e^{-15t}, \quad & t \in (0,T)\\
\vspace{0.1cm}
\displaystyle \frac{\partial u}{\partial x}(L,t) = 0, \quad & t \in (0,T)\\
u(x,0)=0, \; w(x,0)= 0, \quad & x \in (0,L),
\end{cases}
\end{equation}
where $L = 1$, $T = 2$, $\gamma = 2$ and $\beta = 0.5$. The parameter $\mu$ ($n_{\mu} = 1$) belongs to the parameter space $\mathcal{P} = 5 \cdot  [10^{-3}, 10^{-2}]$. This system consists in a parametrized version of the Monodomain equation coupled with the FitzHugh-Nagumo cellular model which describes the excitation-relaxation of the cell membrane in the cardiac tisuue \citep{fitzhugh1961impulses, nagumo1962anactive}. In such a model, the ionic current is a cubic function of the electrical potential $v$ and linear in the recovery variable $w$. Eqs (\ref{1DMonodomain}) have been discretized in space through linear finite elements by considering $N_h = 256$ grid points. We use a one-step, semi-implicit, first order scheme similar to the one discussed in \cite{pagani2018numerical} for time discretization and the treatment of the nonlinear term\footnote{The \texttt{Matlab} library used to compute snapshots and  the numerical results regarding the (local) RB method for problem (\ref{1DMonodomain}) is freely available at \texttt{https://github.com/StefanoPagani/LocalROM}}. The solution of the former problem consists in a parameter-depending travelling wave, which exhibits sharper and sharper fronts as the parameter $\mu$ gets smaller (see \figurename~\ref{fom_solution_1Dmonodomain}).  
We consider $N_{train} = 20$ training-parameter instances uniformly distributed in the parameter space $\mathcal{P}$ and $N_{test} = 19$ testing-parameter instances, each of them corresponding to the midpoint between two consecutive training parameter instances.

\begin{figure}[ht]
\centering
\includegraphics[scale=0.34]{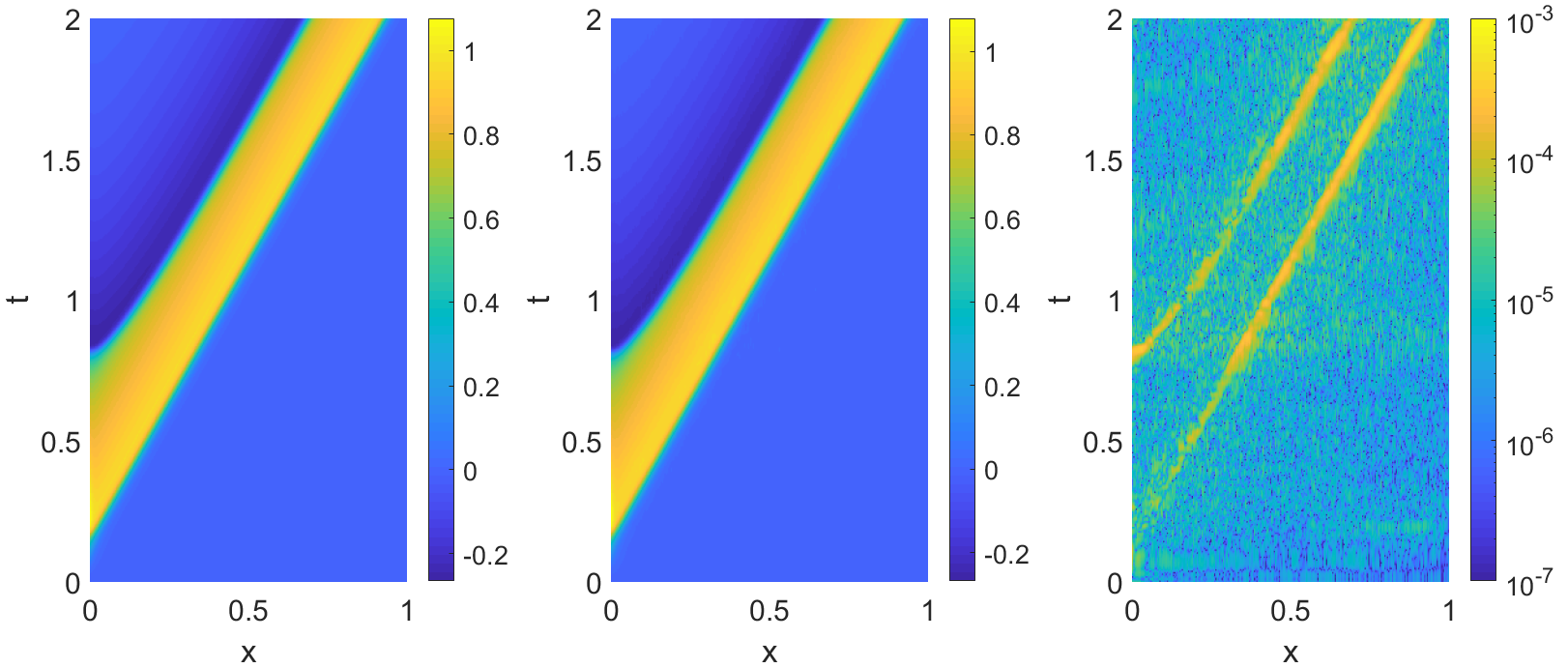}
\caption{\textit{Test 3}: FOM solution (left), DL-ROM solution with $n = 2$ (center) and error indicator $\boldsymbol{\epsilon}_k$ (right) for the testing-parameter instance $\mu_{test} = 0.0062$ in the space-time domain.}
\label{comparison_1D_monodomain}
\end{figure}

\figurename~\ref{comparison_1D_monodomain} shows the FOM solution and the DL-ROM one obtained by setting $n=2$, the dimension of the solution manifold, for the testing-parameter instance $\mu_{test} = 0.0062$. We also report in \figurename~\ref{comparison_1D_monodomain} the error indicator $\boldsymbol{\epsilon}_k$ (\ref{error}), which is higher in correspondence of the large gradients of the FOM solution.

\begin{figure}[ht!]
\centering
\includegraphics[scale=0.15]{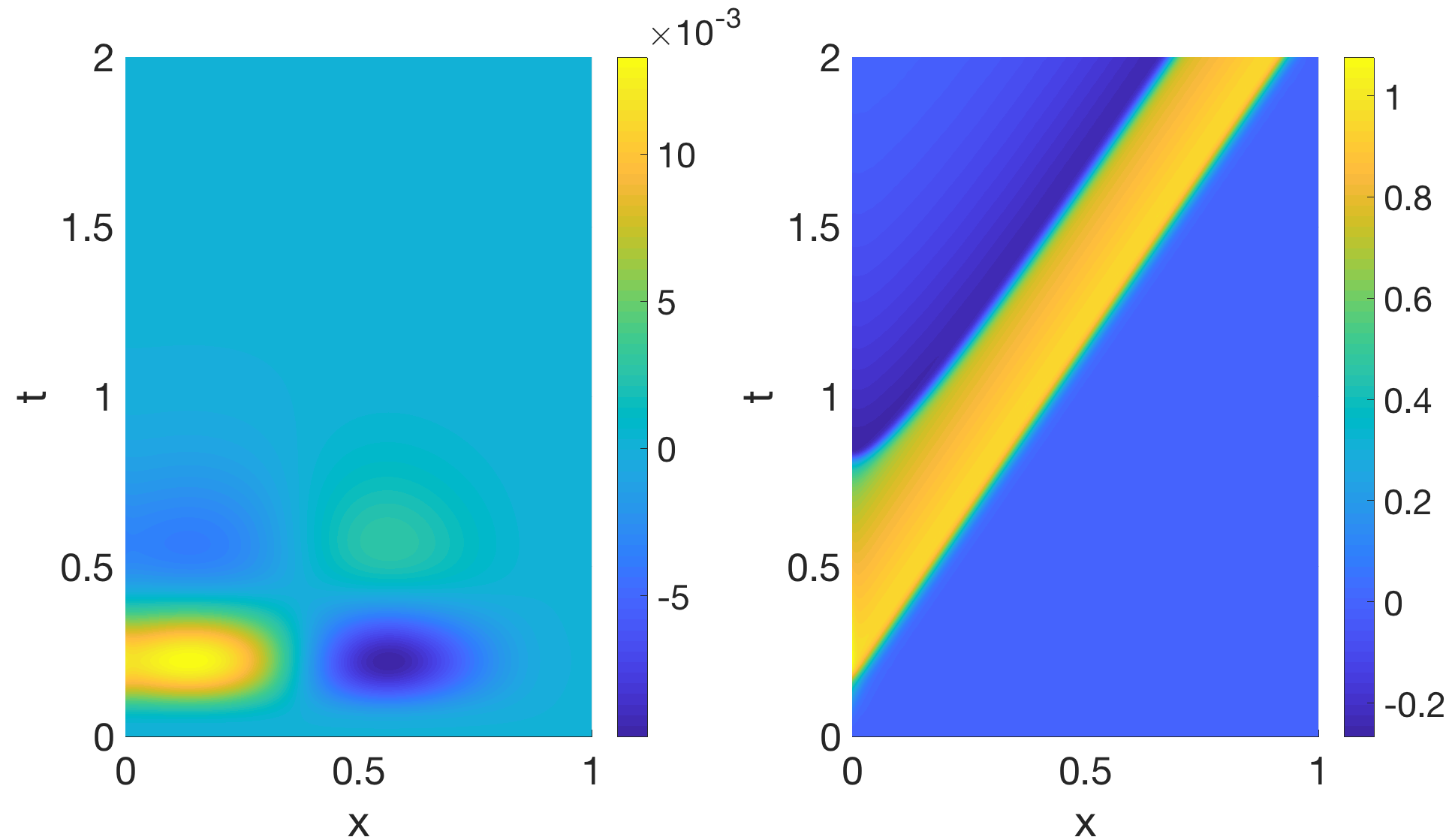}
\caption{\textit{Test 3}: RB solutions for the testing parameter instance $\mu_{test} = 0.0062$ with $n=2$ (left) and $n=66$ (right).}
\label{1D_monodomain_RB}
\end{figure}

\begin{figure}[ht!]
\center{
\includegraphics[scale=0.15]{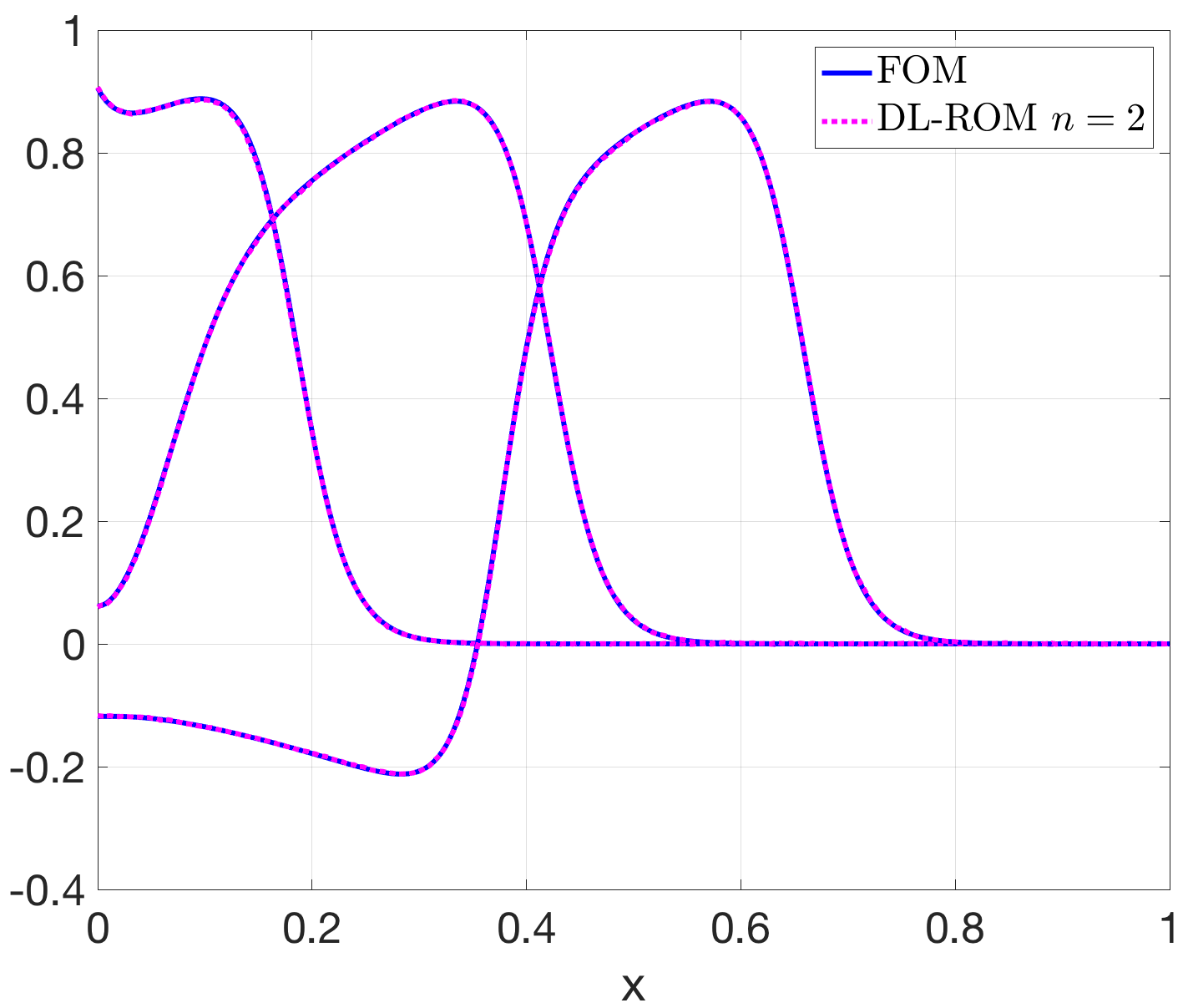}
\hspace{0.3cm}
\includegraphics[scale=0.15]{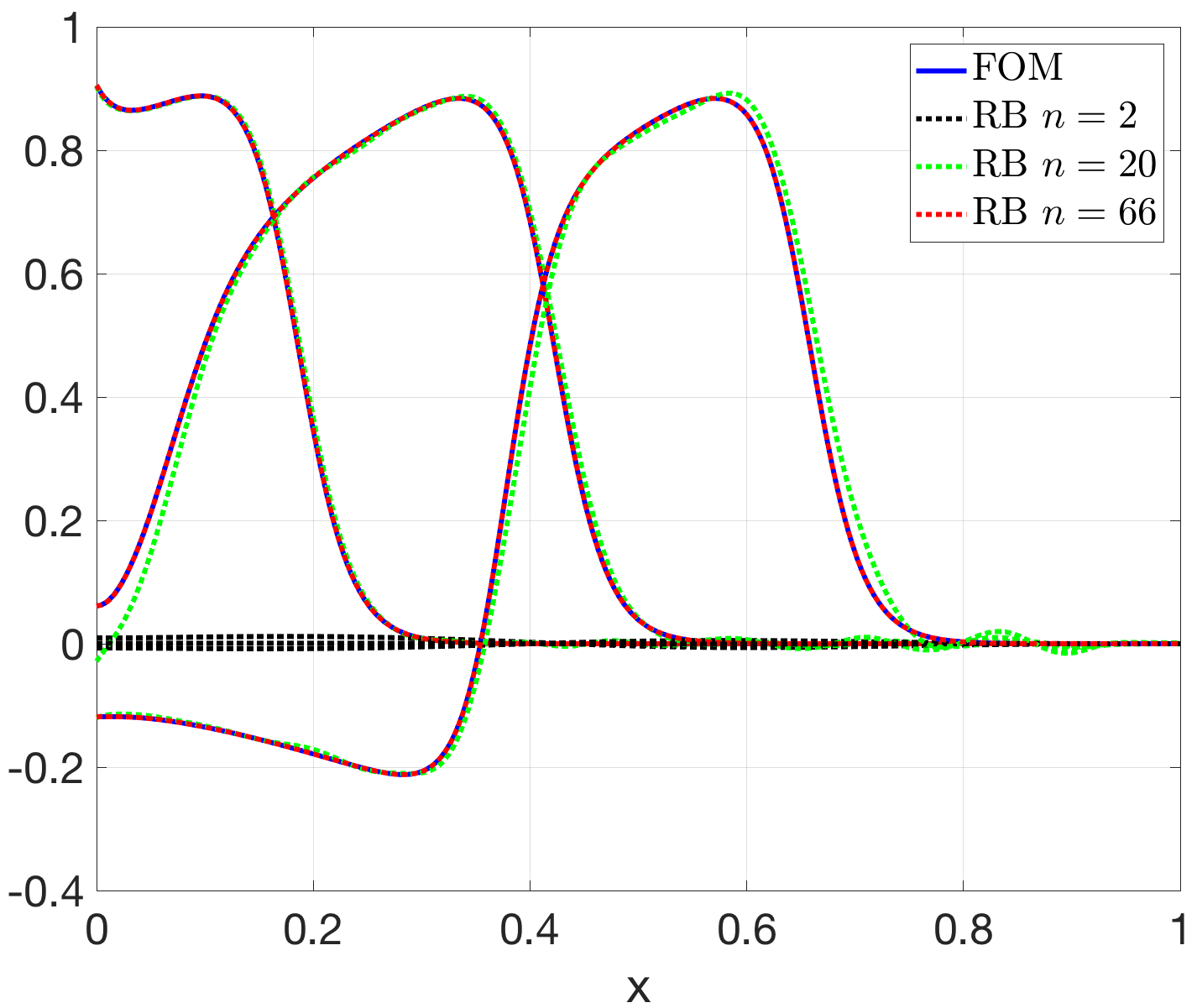}
}
\caption{\textit{Test 3}: FOM and DL-ROM solutions (left) and FOM and RB solutions (right) for the testing-parameter instance $\mu_{test} = 0.0157$ at $t = 0.4962, 0.9975$ and 1.4987.}
\label{comparison_1D_monodomain_time}
\end{figure}

The accuracy obtained by our DL-ROM technique, with $n = 2$, on the testing set is $\epsilon_{rel} = 3.42 \cdot 10^{-3}$. In order to assess the performance of DL-ROM with respect to a linear ROM technique we point out in \tablename~\ref{tab:RB} the maximum number of basis functions among all the clusters, i.e. the dimension of the largest linear trial manifold, required by the (local) RB method in order to achieve the same accuracy obtained through DL-ROM. By increasing the number of clusters, the dimension of the largest linear trial subspace decreases; this does not hold as long as the number of clusters is larger than $k=32$. Indeed, the dimension of some linear subspaces become so small that the error increases with respect the one obtained with fewer clusters. In particular, in \figurename~\ref{1D_monodomain_RB} the RB solutions obtained by considering $n=2$ and $n=66$ basis functions are shown. In \figurename~\ref{comparison_1D_monodomain_time} we compare the FOM solution with the DL-ROM one, obtained for $n = 2$, and the FOM solution with the RB one by setting $n = 2, 20$ and 66, for $\mu_{test} = 0.0157$ at $t = 0.4962, 0.9975$ and 1.4987.
\begin{table}[ht!]
\centering
\begin{tabular}{|c|c|c|c|c|c|}
\hline
k= 1 & k = 2 & k = 4 & k = 8 & k = 16 & k = 32\\
\hline
66 & 68 & 55 & 34 & 26 & 20\\
\hline
\end{tabular}
\caption{\textit{Test 3}: Maximum number of basis functions for the (local) RB method.}
\label{tab:RB}
\end{table}

\begin{figure}[ht!]
\centering
\includegraphics[scale=0.15]{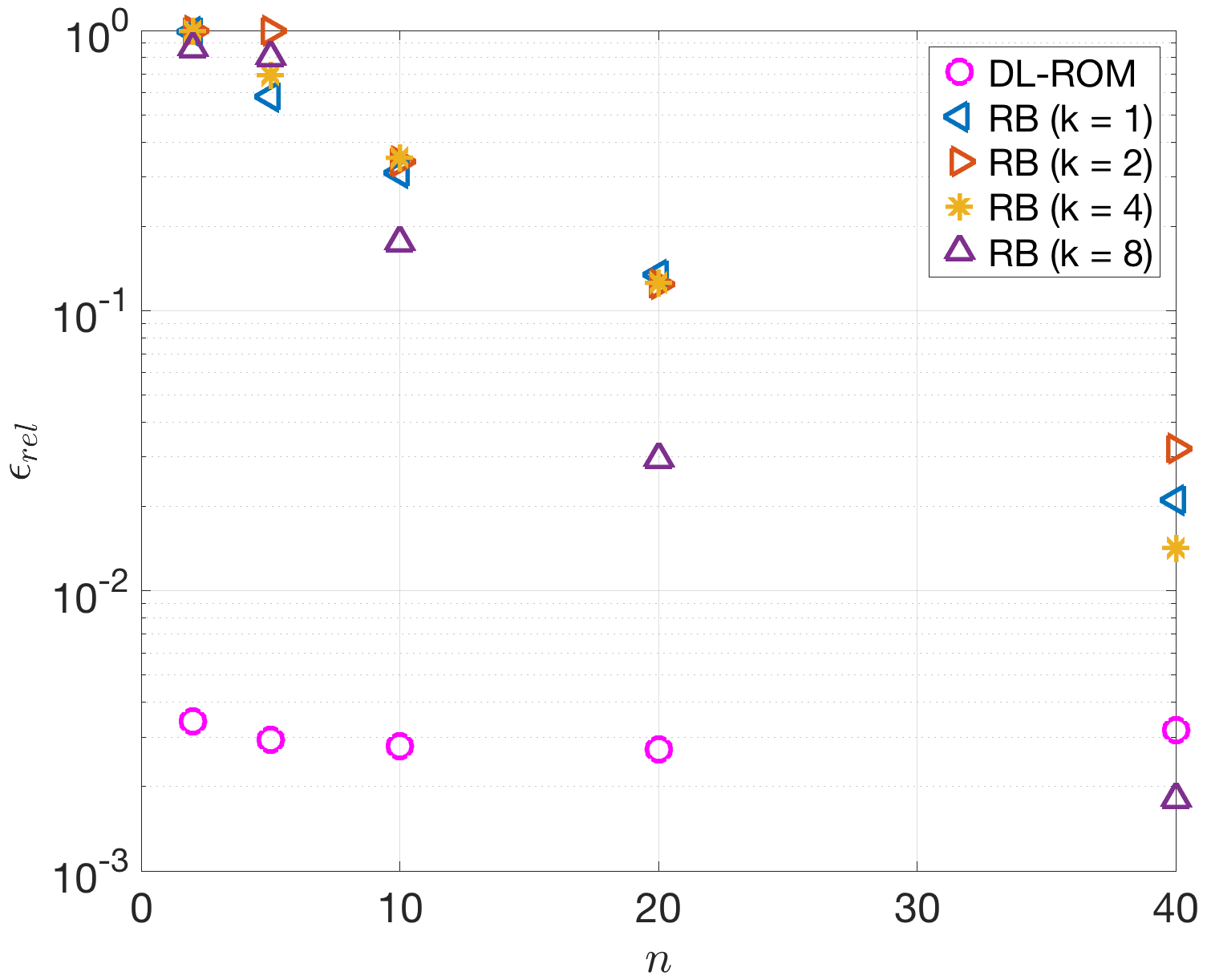}
\caption{\textit{Test 3}: Error indicator $\epsilon_{rel}$ vs. $n$ on the testing set.}
\label{1D_monodomain_convergence}
\end{figure}

The convergence of the error indicator (\ref{relative_error}) as a function of the reduced dimension $n$ is shown in \figurename~\ref{1D_monodomain_convergence}. For the (local) RB method, by increasing the dimension of the largest linear trial manifold, the error indicator decreases, this occurs also by applying the DL-ROM technique for $n \le 20$. The decay of the  error indicator in the latter case is not so remarkable for the same reason pointed out  {in Test 2.1}. If we consider larger values of $n$, e.g. $n = 40$, overfitting occurs, meaning that the neural network model is too complex with respect to the amount of data provided it. For this reason, by considering, for example $n = 40$, the error indicator  $\epsilon_{rel}$ increases.

Finally, in \figurename~\ref{1D_monodomain_convergence_wrt_params} we report the behavior of the loss function and of the error indicator (\ref{relative_error}) with respect to the number of training-parameter instances, i.e. the size of the training dataset. By providing more data to the DL-ROM neural network, its approximation capability {increases, thus yielding  a decrease in the generalization error and the error indicator. In particular, the decay of the loss function with respect to the number of training-parameter instances $N_{train}$ is approximately proportional to $1/N_{train}^3$ and the one of the error indicator  (\ref{relative_error}) is about $1/N_{train}^2$.
\begin{figure}[ht!]
\centering
\includegraphics[scale=0.15]{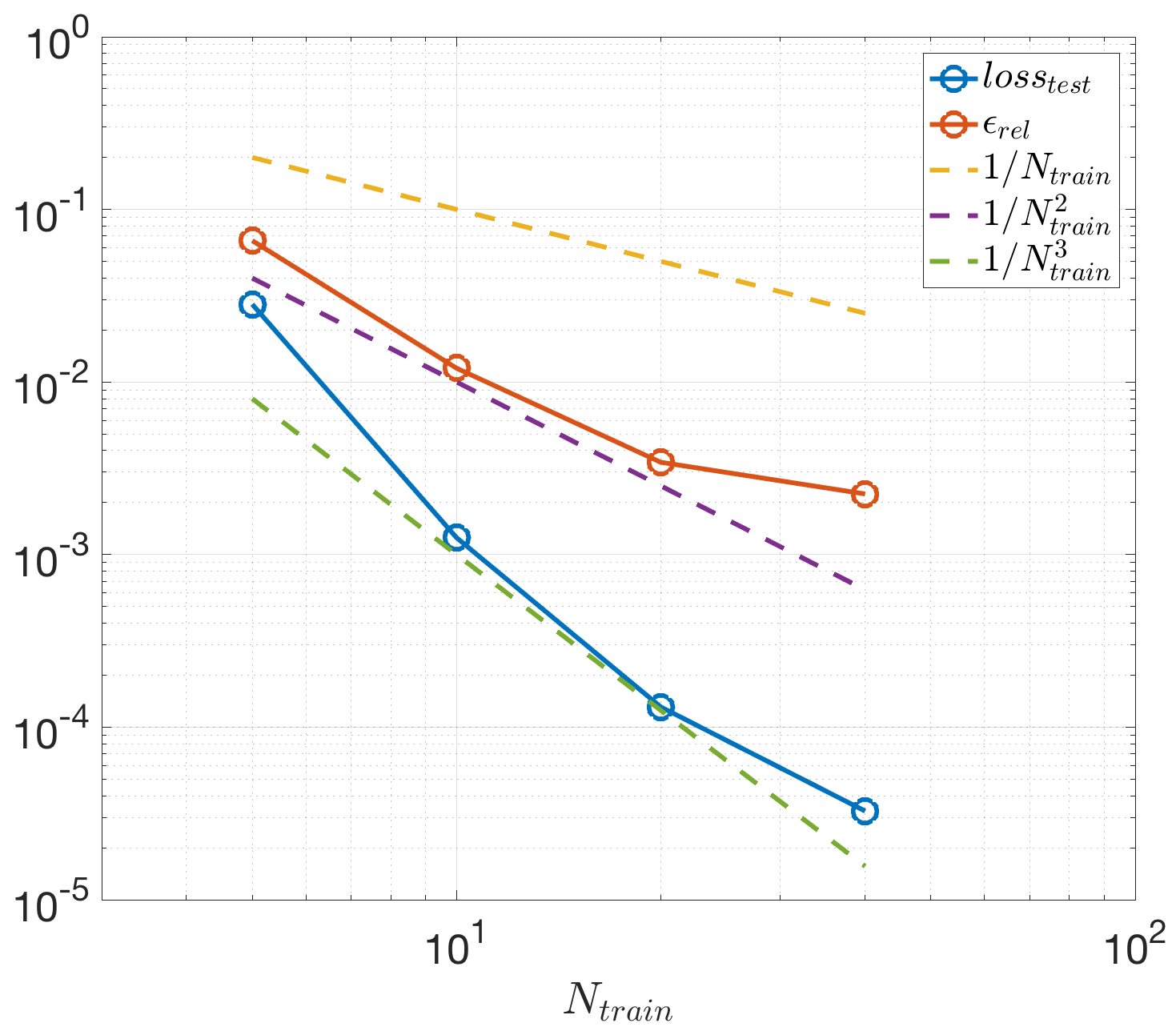}
\caption{\textit{Test 3}: Loss and error indicator $\epsilon_{rel}$ on the testing set vs. number of training-parameter instances of the parameter $\mu$.}
\label{1D_monodomain_convergence_wrt_params}
\end{figure}

\begin{rem}
(\textit{Hyperparameters Tuning}). In order to perform hyperparameters tuning we follow the same procedure used for Test 2.1. We start from the default configuration and we tune the size of the (transposed) convolutional kernels in the (decoder) encoder function, the number of hidden layers in the feedforward neural network and the number of neurons for each hidden layer. In Figure \ref{1D_monodomain_kernel_size} we show the impact of the different hyperparameters on the validation and testing losses. The final configuration of the DL-ROM neural network is the one provided in \tablename~\ref{1D_monodomain_final_configuration}.
\begin{figure}[ht!]
\centering
\includegraphics[width=0.32\textwidth]{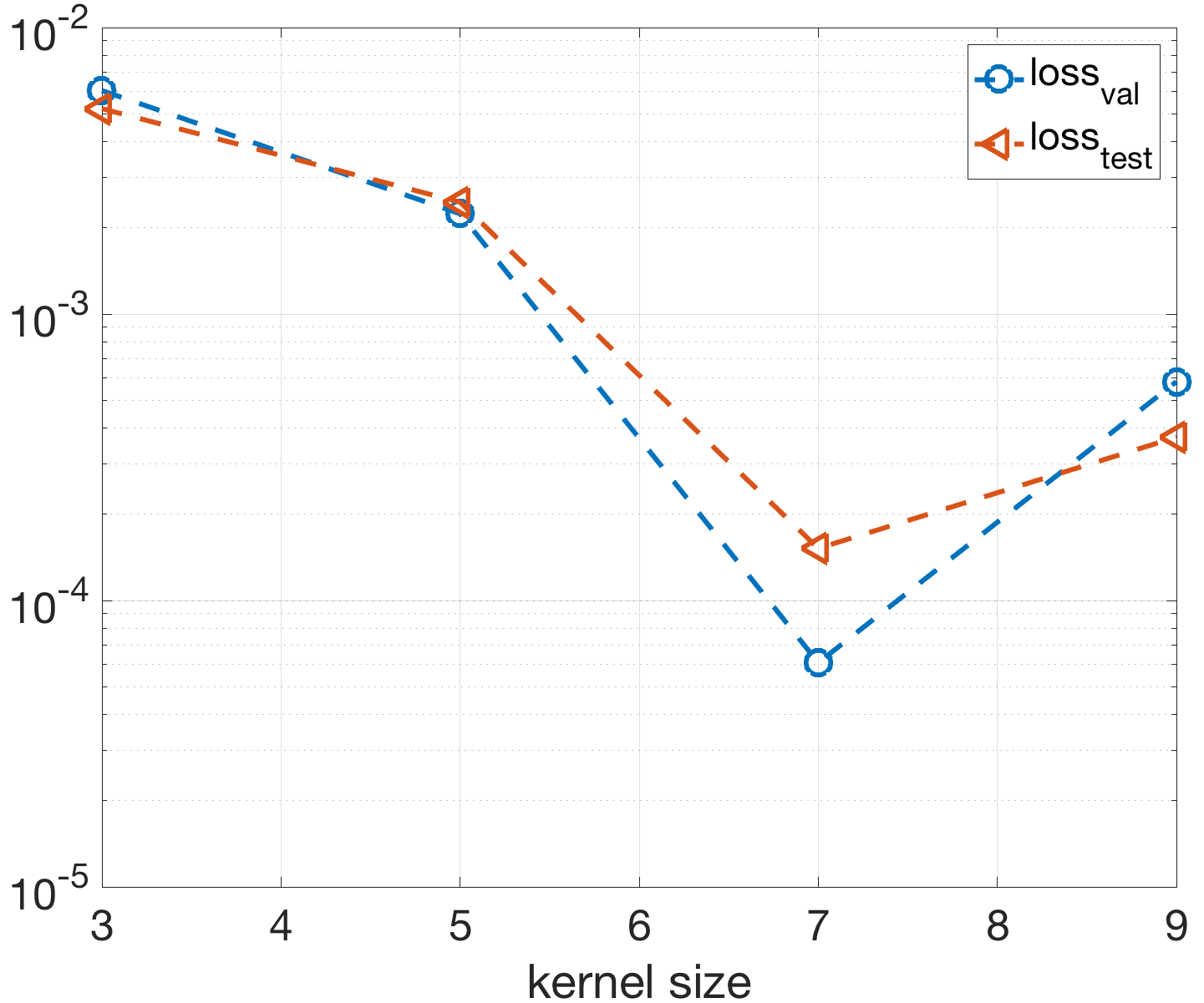}
\includegraphics[width=0.32\textwidth]{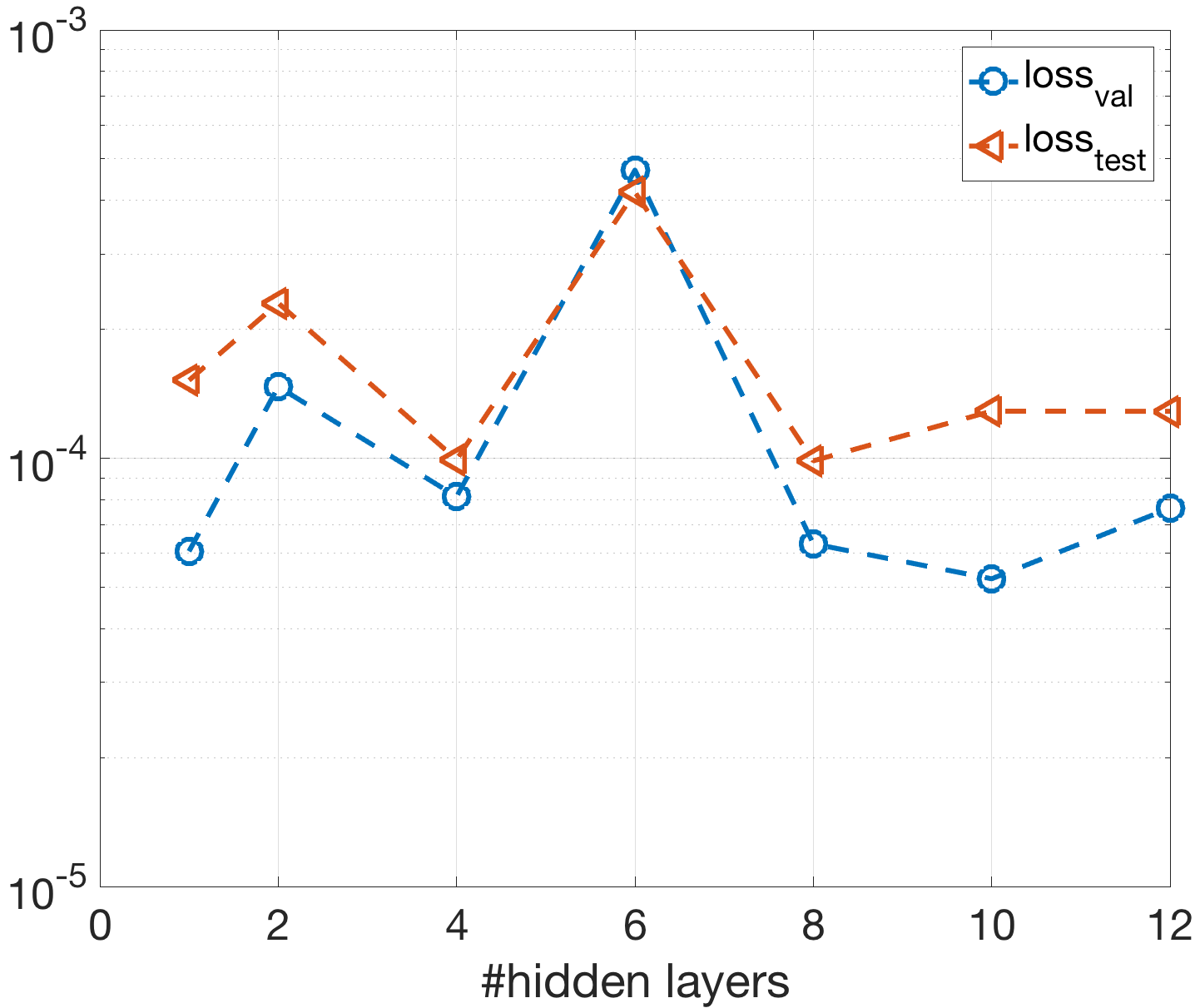}
\includegraphics[width=0.32\textwidth]{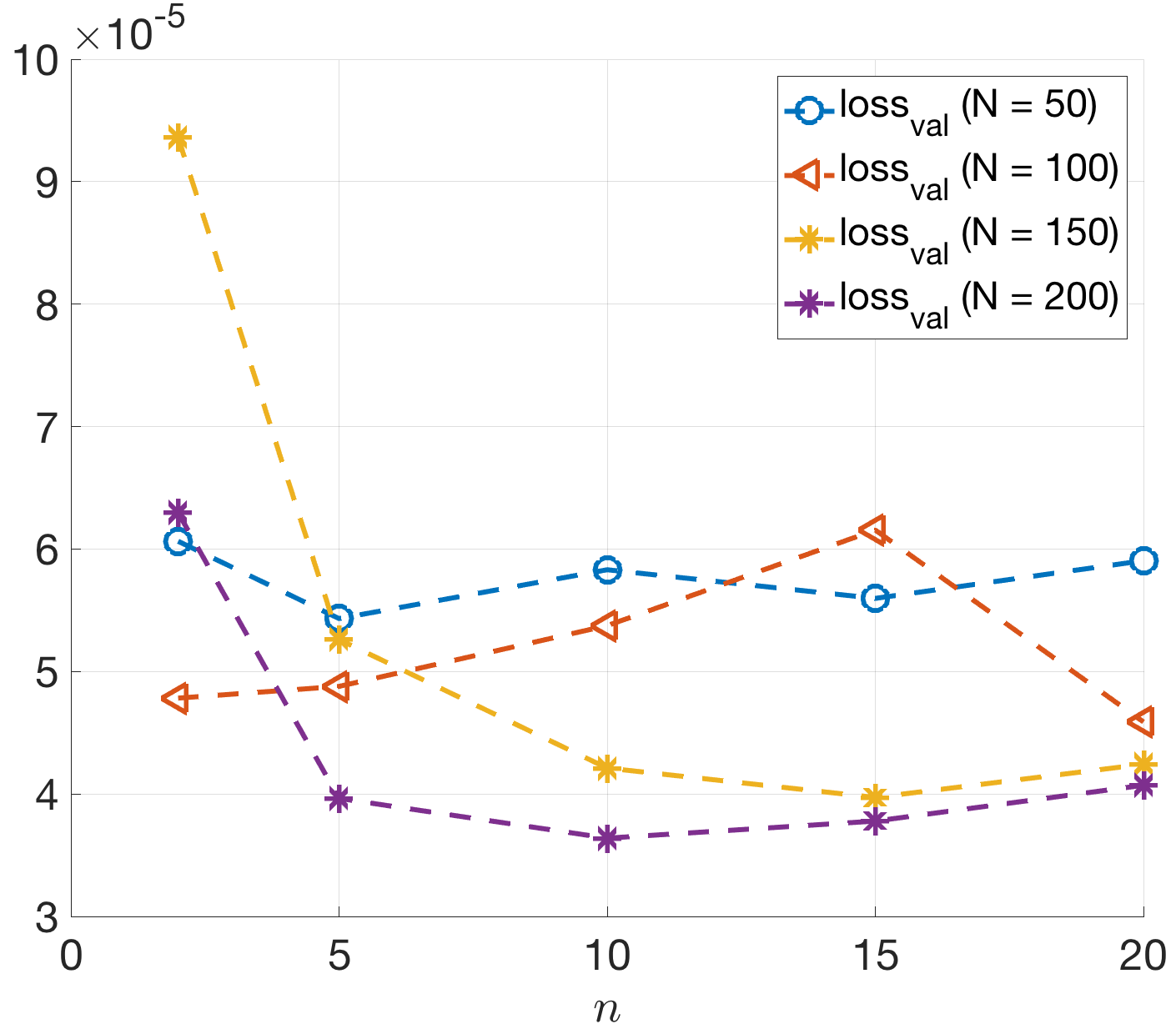}
\caption{\textit{Test 3}: Impact of the kernel size (left), the number of hidden layers (center) and the number of neurons (right) on the validation and testing loss.}
\label{1D_monodomain_kernel_size}
\end{figure}
\begin{table}[ht!]
\begin{center}
\begin{tabular}{|c|c|c|}
\hline
Kernel Size & $\#$ Hidden Layers & $\#$ Neurons \\
\hline
[7, 7] & 1 & 200 \\
\hline
\end{tabular}
\end{center}
\caption{\textit{Test 3}: Final configuration of DL-ROM.}
\label{1D_monodomain_final_configuration}
\end{table}
\end{rem}

\begin{rem}(\textit{Sensitivity with respect to the weight $\omega_h$}). For all the test cases analyzed we set the  parameter $\omega_h$   in the loss function (\ref{loss_encoder})  equal to $\omega_h = 1/2$. In order to justify this choice we performed a sensitivity analysis for problem (\ref{1DMonodomain}) as shown in \figurename~\ref{sensitivity_weights}. For extreme values of $\omega_h$, the error indicator (\ref{relative_error}) worsens of about  one order of magnitude. In particular, not considering the encoder function $\mathbf{f}_n^E$, that corresponds to the case $\omega_h = 1$, yields worse performance of the DL-ROM neural network, as highlighted in \figurename~\ref{sensitivity_weights}. Similarly, by taking $\omega_h = 0$, we would neglect the reconstruction error (that is, the first term in the per-example loss function \eqref{loss_encoder}); this is why the error indicator is large for $\omega_h = 0.1$. All the values of $\omega_h$ in the range $[0.2, 0.9]$ do not yield significant differences  in terms of error indicator, so we decided to set $\omega_h = 1/2$ -- and, as a matter of fact, $ 1-\omega_h = 1/2$.
\begin{figure}[ht]
\centering
\includegraphics[scale=0.14]{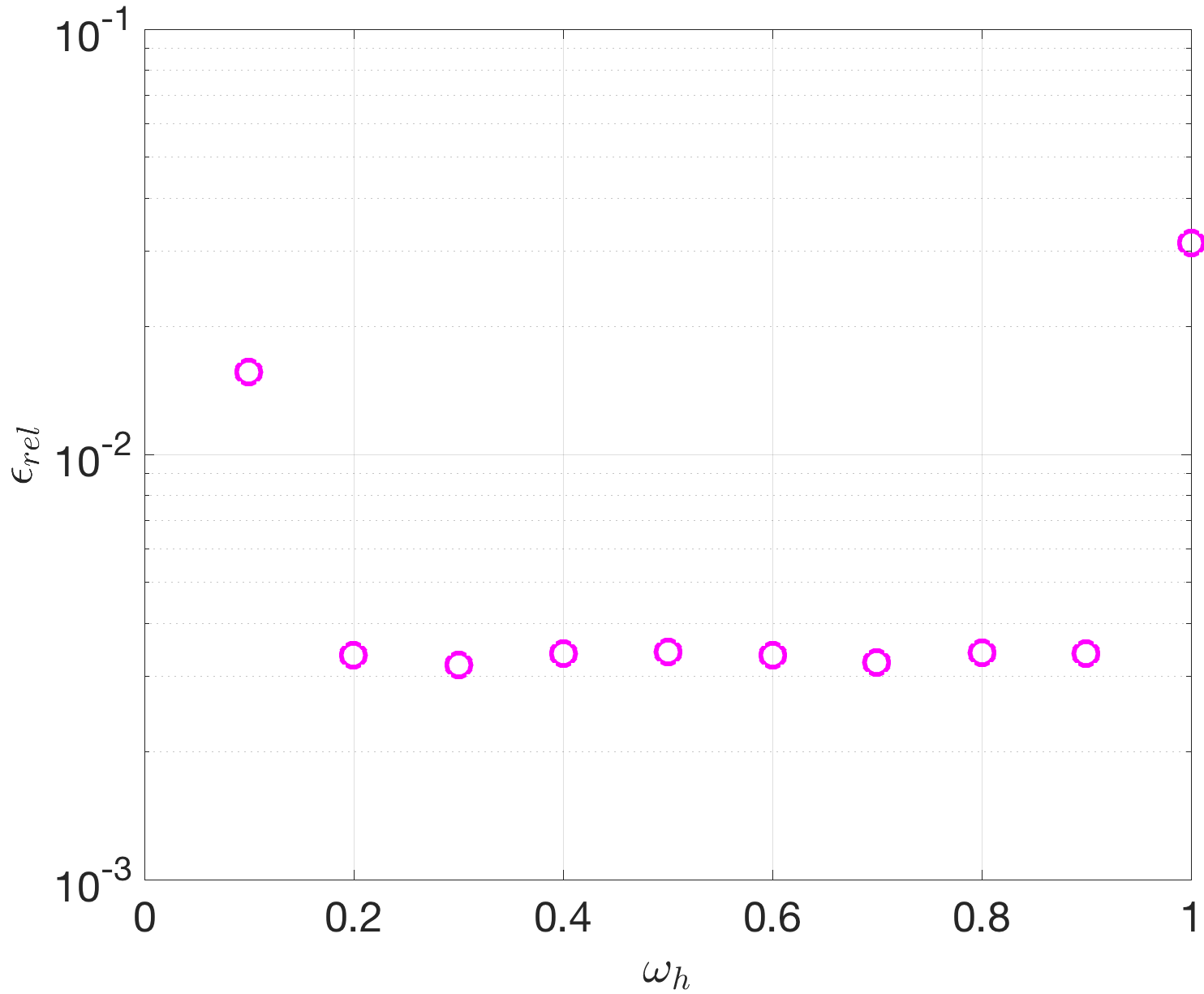}
\caption{\textit{Test 3}: Error indicator $\epsilon_{rel}$ vs. $\omega_h$.}
\label{sensitivity_weights}
\end{figure}
\end{rem}

\section{Conclusions}
\label{sec:5}

In this work we have proposed a novel technique to build low-dimensional ROMs exploiting deep learning models in order to overcome the usual computational bottlenecks shown by classical, linear projection-based ROM techniques (such as the reduced basis method relying on proper orthogonal decomposition) when dealing with problems featuring coherent structures that propagate over time, such as transport and wave-type phenomena, or convection-dominated flows.

The proposed Deep Learning-based Reduced Order Model (DL-ROM) allows to approximate both the solution manifold of a given parametrized nonlinear, time-dependent PDE by means of a low-dimensional, nonlinear trial manifold, and the nonlinear dynamics of the generalized coordinates on such reduced trial manifold, as a function of the time coordinate and the parameters. Both {\em (i)} the nonlinear trial manifold and {\em (ii)} the reduced dynamics are learnt in a non-intrusive way, thus avoiding to query the arrays related to the FOM; the former is learnt by means of the decoder function of a convolutional autoencoder neural network, whereas the latter through a (deep) feedforward neural network, and the encoder function of the convolutional autoencoder.

The numerical results obtained for three different test cases   show that the proposed DL-ROM technique provides sufficiently accurate solutions to the parametrized PDEs involving a low-dimensional solution manifold whose dimension is $n_{\mu}+1$. The proposed DL-ROM outperforms linear ROMs such as the RB method (relying on a global POD basis), as well as  nonlinear approaches exploiting  local POD bases, when  applied both to {\em (i)} problems which are extremely challenging for linear ROMs, such as the linear transport equation or nonlinear diffusion-reaction PDEs coupled to ODEs, and {\em (ii)}  problems which are more tractable using a linear ROM, like Burgers equation, however featuring POD bases with much higher dimension.

Regarding numerical accuracy, the proposed DL-ROM technique provides approximations that are orders of magnitude more accurate than the ones provided by linear ROMs, when keeping the same dimension. We do not obtain remarkable error decays when considering low-dimensional spaces of increasing dimensions, thus making the accuracy of both approximations comparable when dealing with $\mathcal{O}(10^2)$  POD basis functions -- a dimension which makes however linear ROMs infeasible when moving to more involved parametrized problems in higher space dimensions. Regarding computational efficiency, we deem not appropriate to perform comparisons with one-dimensional test cases (on  meshes featuring no more than $\mathcal{O}(10^3)$ degrees of freedom). We will perform the assessment of the computational speedup of our DL-ROM technique compared to linear ROMs in future publications; we expect however to obtain remarkable computational gains when dealing with two and three-dimensional problems for which  linear ROMs are not well-suited to approximate the solution to parametrized, nonlinear time-dependent PDEs. Numerical results shown that DL-ROM allows to generate approximation spaces of dimension close to the intrinsic dimension of the solution manifold, by providing also remarkably improvements in terms of efficiency, will be published in a forthcoming paper.

\section*{Acknowledgments}

We gratefully acknowledge Prof. A. Quarteroni (MOX, Politecnico di Milano) for his stimulating discussions, Dr. S. Pagani (MOX, Politecnico di Milano) for his useful remarks and M. Salvador (MOX, Politecnico di Milano) for kindly sharing the code implementing the FOM of Test 1.

\appendix
\section{Basic concepts of deep learning}
\label{sec:A}

Deep learning (DL) techniques have gained great attention in recent years in several areas like computer vision \citep{krizhevsky2012imagenet, antipov2017face}, natural language processing \citep{sutskever2014sequence, devlin2018bert} and speech recognition \citep{bourlard1989speech, chung2017lip}, due to their ability to discover pattern and extract features from massive datasets, in order to make predictions without providing hand-crafted features. In this section we provide an overview of those deep-learning models which the proposed DL-ROM technique presented in this work relies on.

\subsection{Deep feedforward neural network}

A remarkable example of DL model is the deep feedforward neural network (DFNN). A DFNN is a mathematical function modeling the relationship between a set of input values and some output values \citep{goodfellow2016deep}. This mathematical function is obtained through composition of simpler (nonlinear) functions, or layers, and allows to learn complex hierarchies of features. 
More formally, provided an input $\mathbf{x} \in \mathbb{R}^{N_0}$ a DFNN with $L$ layers takes the form
\begin{equation}
\label{DFNN}
\boldsymbol{\phi}^{DF} : (\mathbf{x}; \boldsymbol{\theta}_{DF}) \mapsto \boldsymbol{\phi}_L(\cdot; \boldsymbol{\theta}_L) \circ \boldsymbol{\phi}_{L-1}(\cdot; \boldsymbol{\theta}_{L-1}) \circ \ldots \circ \boldsymbol{\phi}_1( \mathbf{x}; \boldsymbol{\theta}_1),
\end{equation}
where $\boldsymbol{\phi}_i(\cdot; \boldsymbol{\theta}_i) : \mathbb{R}^{N_i-1} \mapsto \mathbb{R}^{N_i}$, $i = 1, \ldots, L$, refers to the activation function applied at layer $i$ of the DFNN and $\boldsymbol{\theta}_i = (W_i, \mathbf{b}_i)$, with $W_i \in \mathbb{R}^{N_i \times N_{i-1}}$ and $\mathbf{b}_i \in \mathbb{R}^{N_i}$, $i = 1, \ldots, L$, are the weights and the bias of layer $i$ such that $\boldsymbol{\theta}_{DF} = (\boldsymbol{\theta}_1, \ldots, \boldsymbol{\theta}_L)$. We usually refer to the collection of all weights and biases as to the \textit{parameters vector}. Each layer of the network corresponds to a matrix whose values are computed by applying a linear transformation to the previous layer followed by the application of a nonlinear activation function. In particular, referring to \figurename~\ref{feedforward_neural_network}, $\mathbf{y}_0 = \mathbf{x} \in \mathbb{R}^{N_0}$ is the input layer, $\mathbf{y}_L = \boldsymbol{\phi}^{DF}(\mathbf{x}; \boldsymbol{\theta}_{DF}) \in \mathbb{R}^{N_L}$ is the output layer, and each hidden layer $\mathbf{y}_i \in \mathbb{R}^{N_i}$, $i = 1, \ldots, L - 1$, takes the form
\begin{equation*}
\mathbf{y}_i = \boldsymbol{\phi}_i(W_i \mathbf{y}_{i-1} + \mathbf{b}_i).
\end{equation*}
\begin{figure}[ht!]
\centering
\includegraphics[scale=0.3]{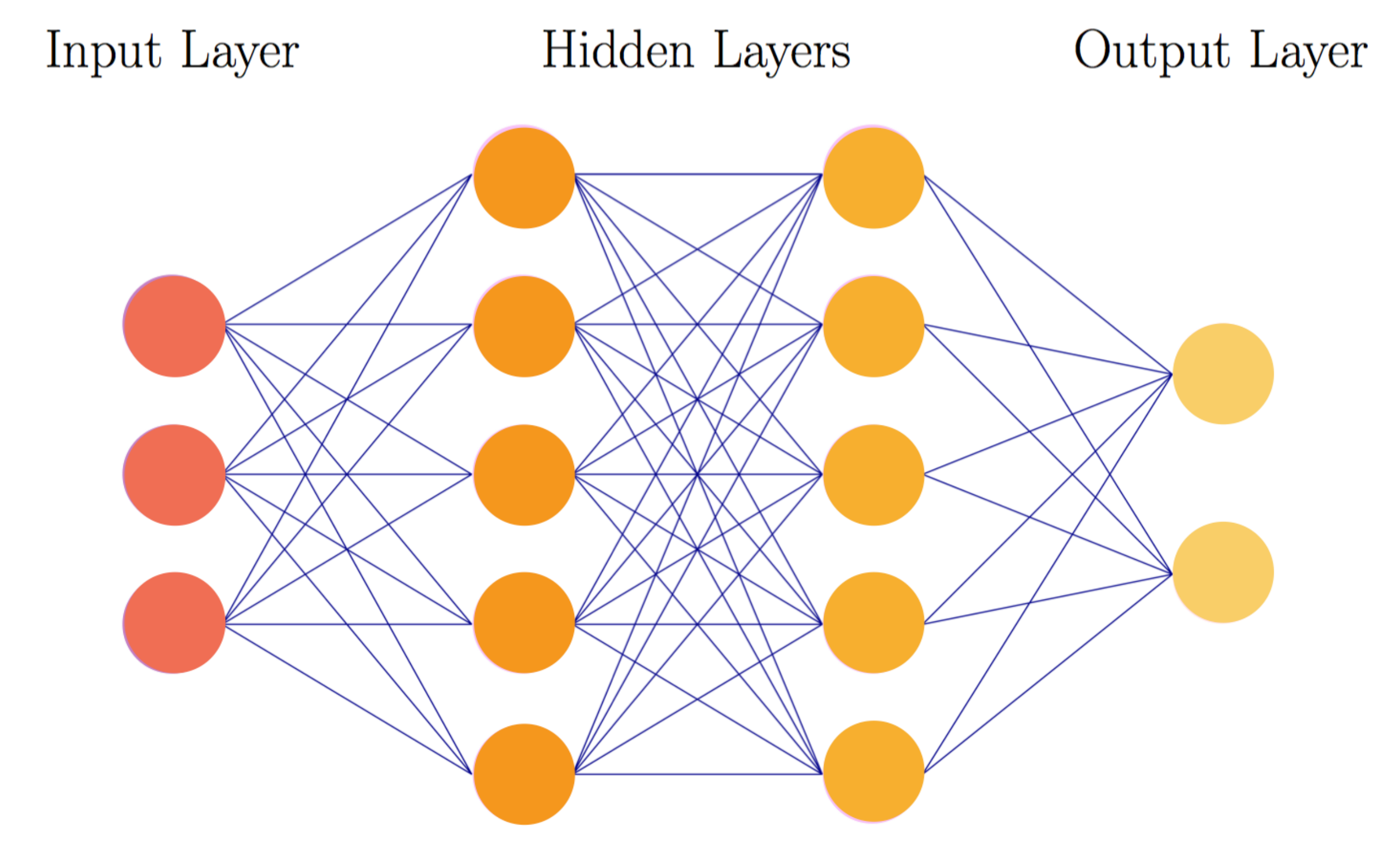}
\caption{Feedforward neural network}
\label{feedforward_neural_network}
\end{figure}

Given a set of $M$ input-output pair observations $\{( \mathbf{x}^i, \mathbf{y}^i) \}_{i=1}^M$ and considering a \textit{supervised learning} paradigm \citep{goodfellow2016deep}, the learning task consists in finding the optimal parameters vector $\boldsymbol{\theta}_{DF}^*$ by solving the optimization problem 
\begin{equation}
\label{minimization_problem}
\min_{\boldsymbol{\theta}_{DF}} \mathcal{J}(\boldsymbol{\theta}_{DF}) = \min_{\boldsymbol{\theta}_{DF}} \frac{1}{M}\sum_{i=1}^M \mathcal{L}(\mathbf{y}^i, \mathbf{y}_L^i; \boldsymbol{\theta}_{DF})
\end{equation}
where $\mathcal{J}$ is the loss (or cost) function, and $\mathcal{L}$ is the per-example loss function, measuring the mismatch between the desired observed output $\mathbf{y}^i$ and the approximated one $\mathbf{y}_L^i$. Problem (\ref{minimization_problem}) is usually solved by means of the gradient descent method exploiting the back-propagation algorithm \citep{rumelhart2986learning} to compute the derivatives of the loss function with respect to parameters. In particular, the gradient descent method requires to evaluate 
\begin{equation}
\nabla_{\boldsymbol{\theta}_{DF}} \mathcal{J}(\boldsymbol{\theta}_{DF}) = \frac{1}{M}\sum_{i=1}^M \nabla_{\boldsymbol{\theta}_{DF}}\mathcal{L}(\mathbf{y}^i, \mathbf{y}_L^i; \boldsymbol{\theta}_{DF}), 
\label{gradient_loss}
\end{equation}
a task which might easily become prohibitive when the size $M$ of the training dataset is very large, thus causing a single step of the gradient descent method to require a huge amount of time. The stochastic gradient descent (SGD) method allows to reduce  the computational cost associated to the computation of the gradient of the loss function, by exploiting the fact that (\ref{gradient_loss}) can be considered as an expectation over the entire training dataset. Such an expectation can be approximated using a small set (or \textit{minibatch}) of samples; hence,   at each iteration the SGD method samples a minibatch of $m < M$ data points, 
 drawn (e.g., uniformly) from the training dataset \citep{goodfellow2016deep}, and approximates  the gradient \eqref{gradient_loss} of the loss function  by
\begin{equation*}
\widehat{\nabla}_{\boldsymbol{\theta}_{DF}} \mathcal{J}(\boldsymbol{\theta}_{DF}) = \frac{1}{m}\sum_{i=1}^m \nabla_{\boldsymbol{\theta}_{DF}}\mathcal{L}(\mathbf{y}^i, \mathbf{y}_L^i; \boldsymbol{\theta}_{DF}). 
\end{equation*}

\subsection{Convolutional neural network}
\label{Convolutional Neural Network}

Convolutional neural networks (CNNs) \citep{lecun1998gradient} are the standard neural network architecture in computer vision tasks, since they are well-suited to high-dimensional and spatially distributed data like images. This is due to the local approach of convolutional layers which enables them to exploit spatial correlations among pixels in order to extract low-level features of the input to carry out the task. The main ingredients of a convolutional layer are convolutional kernels, or filters, which consist in tensors of smaller dimensions with respect to the input. Each element of a feature map is obtained by sliding the kernel over the image and by computing the discrete convolution, as shown in \figurename~\ref{convolutional_neural_network}. \\
\begin{figure}[ht!]
\centering
\includegraphics[scale=0.275]{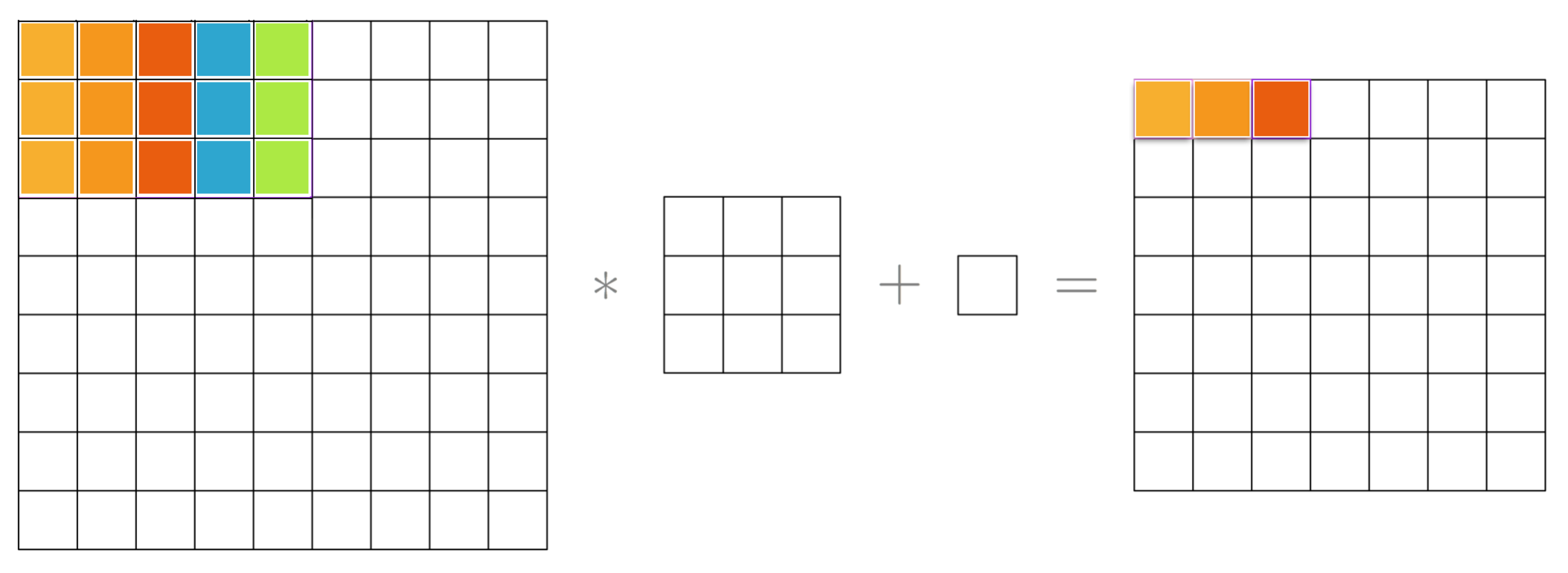}
\caption{Computation of the elements of a feature map in a convolutional layer.}
\label{convolutional_neural_network}
\end{figure}

Considering a 3-dimensional input $Y_0 = X \in \mathbb{R}^{N_0^1 \times N_0^2 \times N_0^3}$ and a bank of $K_i$ convolutional filters in layer $i$ denoted as $W_i^k \in \mathbb{R}^{n_i^1 \times n_i^2 \times n_i^3}$, $i = 1, \ldots, L$ and $k = 1, \ldots, K_i$, the $k$-th feature map is computed as
\begin{equation*}
Y_i^k = \boldsymbol{\phi}_i(W_i^k \ast Y_{i-1} + b_i^k).
\end{equation*}
where $Y_i \in \mathbb{R}^{N_i^1 \times N_i^2 \times N_i^3}$ (or, equivalently, $Y_i^k \in \mathbb{R}^{N_i^1 \times N_i^2}$) with $N_i^1$ and $N_i^2$ depending on $n_i^1$ and $n_i^2$, respectively, the padding and the striding strategies, and $N_i^3 = K_i$.

Convolutional layers are characterized by \textit{shared parameters}, that is, weights are shared by all the elements (neurons) in a particular feature map, and   \textit{local connectivity}, that is, each neuron in a feature map is connected only to a local region of the input. Parameter sharing allows convolutional layers to enjoy another property: translation invariance or, more precisely, translation equivariance. This means that if the input varies, the output changes accordingly \cite{goodfellow2016deep}. In particular, if we apply a transformation to the input $Y_0$ and then compute the convolution, the result is the same we would obtain by computing the convolution and then applying the transformation to the output. The two properties above increase efficiency of CNNs, both in terms of memory and computational costs, with respect to DFNNs, thus making them preferable to the latter when dealing with extremely high-dimensional data.

\subsection{Autoencoder neural network}

Autoencoders (AEs) \citep{bourlard1998auto-association, hinton1994autoencoders} are a particular type of feedforward neural networks aiming at learning, under suitable constraints, the identity function
\begin{equation}
\label{autoencoder}
\mathbf{f}^{AE} (\cdot; \boldsymbol{\theta}_E, \boldsymbol{\theta}_D) : \mathbf{x}_h \mapsto \mathbf{\tilde{x}}_h \quad \textnormal{with} \quad \mathbf{\tilde{x}}_h \simeq \mathbf{x}_h.
\end{equation}
Internally, an autoencoder has a hidden layer consisting in a \textit{code} used to represent the input. We focus on undercomplete autoencoders \citep{goodfellow2016deep} where the constraint imposed is the reduction of the dimension of the code with respect to the input and output dimension.

By considering the input $\mathbf{y}_0 = \mathbf{x}_h \in \mathbb{R}^{N_h}$ and the output $\mathbf{y}_L = \mathbf{\tilde{x}}_h \in \mathbb{R}^{N_h}$, an autoencoder is composed by two main parts (see \figurename~\ref{autoencoder_neural_network})
\begin{itemize}
\item the \textit{encoder} function $\mathbf{f}_{n}^E(\cdot; \boldsymbol{\theta}_E) : \mathbf{x}_h \mapsto \mathbf{\tilde{x}}_n = \mathbf{f}_{n}^E( \mathbf{x}_h; \boldsymbol{\theta}_E)$, where $\mathbf{f}_{n}^E (\cdot; \boldsymbol{\theta}_E) : \mathbb{R}^{N_h} \rightarrow \mathbb{R}^{n}$ and $n \ll N_h$, mapping the high-dimensional input $\mathbf{x}_h$ onto the low-dimensional code $\mathbf{\tilde{x}}_n$. The encoder function depends on a vector of parameters $\boldsymbol{\theta}_E \in \mathbb{R}^{N_E}$ collecting all the weights and biases specifying the function itself;

\item the \textit{decoder} function $\mathbf{f}_{h}^D(\cdot; \boldsymbol{\theta}_D) : \mathbf{\tilde{x}}_n \mapsto \mathbf{\tilde{x}}_h = \mathbf{f}_{h}^D(\mathbf{\tilde{x}}_n; \boldsymbol{\theta}_D)$, where $\mathbf{f}_{h}^D(\cdot; \boldsymbol{\theta}_D) : \mathbb{R}^{n} \rightarrow \mathbb{R}^{N_h}$, mapping the code $\tilde{\mathbf{x}}_n$ to an approximation of the original high-dimensional input $\mathbf{\tilde{x}}_h$. Similarly to the encoder function, the decoder function depends on a vector of parameters $\boldsymbol{\theta}_D \in \mathbb{R}^{N_D}$ collecting all the weights and biases specifying the function itself.
\end{itemize}
The autoencoder is then defined as $$\mathbf{f}^{AE} (\cdot; \boldsymbol{\theta}_E, \boldsymbol{\theta}_D) : \mathbf{x}_h \mapsto  \tilde{\mathbf{x}}_h = \mathbf{f}_{h}^D( \mathbf{f}_{n}^E( \mathbf{x}_h; \boldsymbol{\theta}_E); \boldsymbol{\theta}_D).$$ 
\begin{figure}[ht!]
\centering
\includegraphics[scale=0.3]{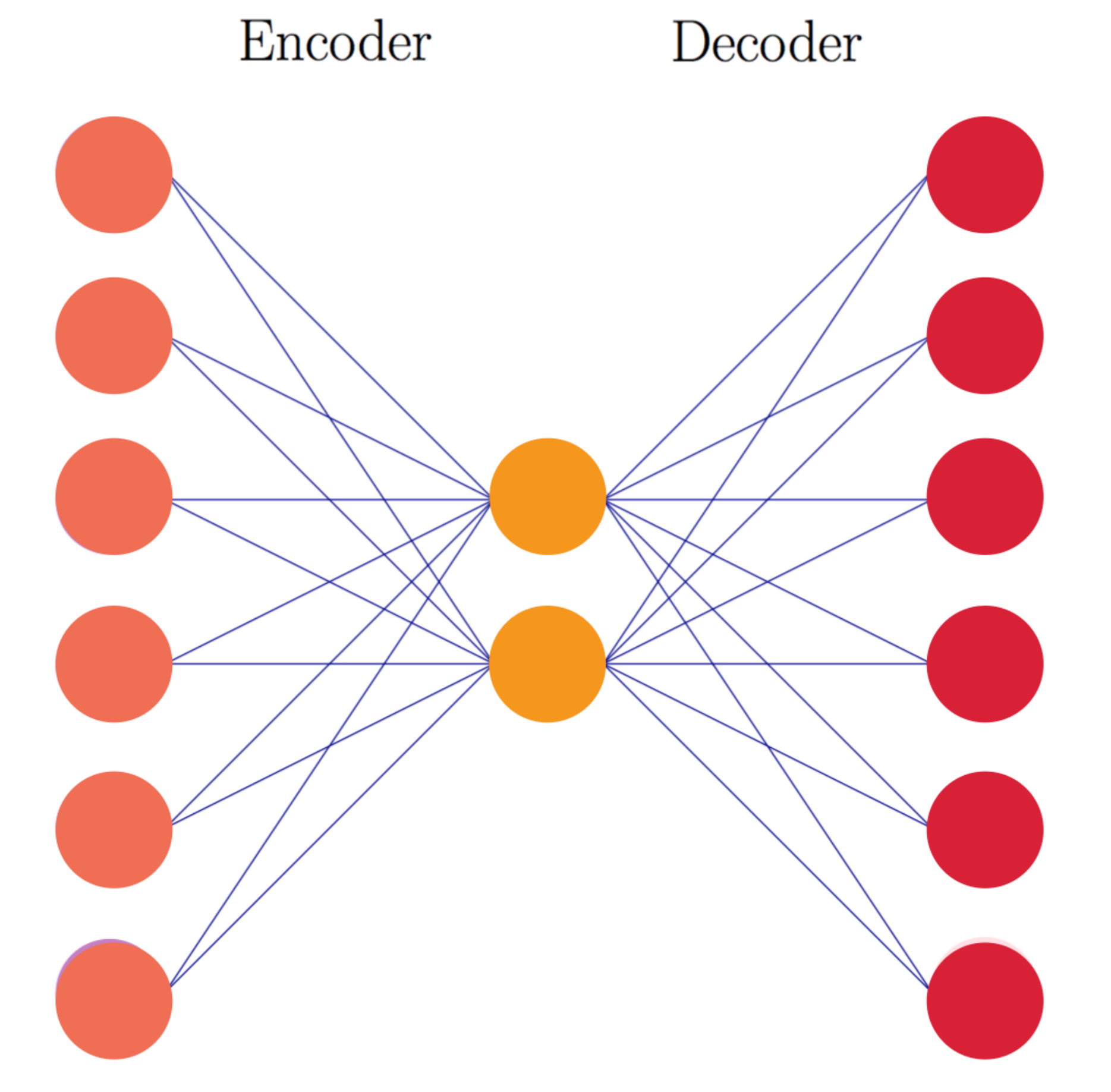}
\caption{Autoencoder neural network.}
\label{autoencoder_neural_network}
\end{figure}

Autoencoder learning lays within the \textit{unsupervised learning} paradigm \citep{goodfellow2016deep}  since its goal is to reconstruct the input being the target output an approximation of the input. An autoencoder not only learns a low-dimensional representation of the high-dimensional input but also learns how to reconstruct the input from the code through the encoder and the decoder functions. 

When dealing with large inputs, as the ones arising from the discretization of system (\ref{FOM}), the use of a feedforward autoencoder may become prohibitive as the number of parameters (weights and biases) required may be very large. As pointed out in \ref{Convolutional Neural Network}, parameter sharing and local connectivity allow to reduce the numbers of parameters of the network and the number of associated computations, both in the forward and in the backward pass, hence the idea of relying on convolutional autoencoders for the sake of building our DL-ROM technique.


\bibliography{Bibliography}

\end{document}